\documentclass[12pt]{amsart}
\usepackage{latexsym, amssymb,enumerate}
\usepackage{appendix}

\newcommand{\be}{\begin{equation}}
\newcommand{\ee}{\end{equation}}
\newcommand{\beq}{\begin{eqnarray}}
\newcommand{\eeq}{\end{eqnarray}}
\newcommand{\cM}{\mathcal{M}}
\newcommand{\ud}{\mathrm{d}}
\newcommand{\re}{\mathrm{e}}

\newcommand{\pr}{\partial}
\newcommand{\fr}{\frac}

\newtheorem{prop}{Proposition}[section]
\newtheorem{theo}[prop]{Theorem}
\newtheorem{lemm}[prop]{Lemma}

\newtheorem{rema}[prop]{Remark}

\newtheorem{defi}[prop]{Definition}

\def\begeq{\begin{equation}}
\def\endeq{\end{equation}}

\numberwithin{equation} {section}

\begin{document}

\title {Normalized Ricci flows and conformally compact Einstein metrics}
\begin{abstract}
In this paper, we investigate the behavior of the normalized Ricci flow
on asymptotically hyperbolic manifolds. We show that the normalized Ricci flow exists globally and converges to an Einstein metric when starting from a non-degenerate and sufficiently Ricci pinched metric. More importantly we use maximum principles to establish the regularity of conformal compactness along the normalized Ricci flow including that of the limit metric at time infinity. Therefore we are able to recover the existence results  in \cite{GL} \cite{Lee} \cite{Bi} of conformally compact Einstein metrics with conformal infinities which are perturbations of that of  given non-degenerate conformally compact Einstein metrics.  
\end{abstract}

\renewcommand{\subjclassname}{\textup{2000} Mathematics Subject Classification}

\author { Jie Qing $^\dag$,
  Yuguang Shi$^\ddag$ and Jie Wu $^\ddag$}

\address{Jie Qing, Department of Mathematics, University of California, Santa Cruz, CA 95064, USA} \email{qing@ucsc.edu}

\address{Yuguang Shi, Key Laboratory of Pure and Applied mathematics, School of Mathematics Science, Peking University,
Beijing, 100871, P.R. China.} \email{ygshi@math.pku.edu.cn}

\address{Jie Wu, Key Laboratory of Pure and Applied mathematics, School of Mathematics Science, Peking University,
Beijing, 100871, P.R. China.} \email{wujie@math.pku.edu.cn}

\thanks{$^\dag$ Research partially supported by NSF DMS-1005295 and CNSF 10728103.}
\thanks{$^\ddag$ Research partially supported by  NSF grant of China 10725101 and 10990013.}

\date{May 2011}
\maketitle

\markboth{ Jie Qing  and Yuguang Shi, Jie Wu}{}

\section {introduction}

Since the seminal work of Fefferman and Graham \cite{FG} there have been great interests in the study of conformally compact Einstein metrics. Lately the use of conformally compact Einstein manifolds in the so-called AdS/CFT correspondence in string theory proposed as a promising quantum theory of gravity have accelerated developments of the study of conformally compact Einstein manifolds. As it was foreseen in \cite{FG}, the study of conformally compact Einstein manifolds now becomes one of the most active research area in conformal geometry. But the existence of conformally compact Einstein metrics remains to be a challenging open problem in large. 

In this paper we study the normalized Ricci flows on asymptotically hyperbolic manifolds and use normalized Ricci flows to construct conformally compact Einstein metrics. We recall that Ricci flow starting from a metric $g_0$ on a manifold $\cM^n$ is a family of metrics $g(t)$ that  satisfies the following:
\begin{equation}\label{RF}
\left\{\aligned \frac d{dt} g(t) & = -2  Ric_{g(t)}\\
                           g(0) & = g_0\endaligned\right.
\end{equation}
We then consider the normalized Ricci flow as follows:
\begin{equation}\label{NRF}
\left\{\aligned \frac d{dt} g(t) & = -2  (Ric_{g(t)} + (n-1)g(t)) \\
                           g(0) & = g_0\endaligned\right.
\end{equation}
It is easily seen that \eqref{RF} is equivalent to \eqref{NRF}. In fact explicitly
$$
g^N(t)  = e^{-2(n-1)t}g(\frac 1{2(n-1)}(e^{2(n-1)t}-1))
$$
solves \eqref{NRF} if and only if $g(t)$ solves \eqref{RF}. We like to mention the recent nice work 
\cite{Ba} on Ricci flows on asymptotically hyperbolic manifolds.

Naturally one initial step is to study normalized Ricci flows starting from metrics that are close to be Einstein. Such questions on compact manifolds were studied in \cite{Ye}, where it was observed that the normalized Ricci flow exists globally and converges exponentially to an Einstein metric if the initial metric $g_0$ is sufficiently Ricci pinched and is non-degenerate. There are also several works in the non-compact cases. In [LY] the stability of the hyperbolic space under the normalized Ricci flow was established. The stability of the hyperbolic space in \cite{LY} later is improved and extended in \cite{Bam1} \cite{Bam2} \cite{SSS} \cite{Su}.  

To be more precise we say a metric $g$ on a manifold $\cM^n$ is $\epsilon$-Einstein if 
\begin{equation}\label{Ricci pinch}
\|h_g\|\leq \epsilon
\end{equation}
on $\cM^n$, where the Ricci pinching curvature $h_g = Ric_g + (n-1)g$. The non-degeneracy of a metric is defined to be the first $L^2$ eigenvalue of the linearization of the curavture tensor $h$ as follows:
\begin{equation}\label{non-degeratcy}
\lambda = \inf \frac {\int_ \mathcal{M}  \langle\big(\Delta_L+2(n-1)\big)u_{ij}, u_{ij}\rangle}{  \int_\mathcal{M}\|u\|^2} 
\end{equation}
where the infimum is taken among symmetric $2$-tensors $u$ such that 
$$
\int_{\cM}(|\nabla u|^2+ |u|^2)dv < \infty
$$ 
and $\Delta_L$ is Lichnerowicz Laplacian on symmetric 2-tensors.


We first, based on the ideas in \cite{Ye} \cite{LY}, obtain the following global existence and convergence theorem of the normalized Ricci flow on non-compact manifolds analogous to the ones in \cite{Ye}. 

\begin{theo}\label{maintheorem1}
For any $n\geq3$ and positive constants $k_0,k_1, v_0,\lambda_0$, and $\alpha$ there exists $\epsilon>0$ depending only on $n, k_0, k_1, v_0,\lambda_0$, and $\alpha$ such that the normalized Ricci
flow starting from a metric $g_0$ exists for all the time and converges exponentially to
an Einstein metric, provided that 
\begin{enumerate}[(1)]
\item
$\|Rm_{g_0}\|\leq k_0, \|\nabla Rm_{g_0}\|\leq k_1$
\item
the volume bound $\text{vol}(B_{g_0}(x, 1))  \geq v_0$, for all $x\in \cM^n$
\item
$g_0$ is with non-degeneracy $\geq \lambda_0$
\item
$\int_{\cM} \exp( - \alpha d(x, x_0)) dv_{g_0} < \infty$  
\item
$g_0$ is $\epsilon$-Einstein, and
\item
\begin{equation*}
\int_{\cM} \|h_{g_0}\|^2_{g_0}  dv_{g_0} \leq \epsilon,
\end{equation*}
\end{enumerate}
where $d(x, x_0)$ is the distance to a fixed point $x_0\in \cM^n$.
\end{theo}

As observed in \cite{LY}, one may replace the $L^2$ small condition (6) in Theorem \ref{maintheorem1} by some decay of the Ricci pinching curvature $h$, which is better particularly in the context of asymptotically hyperbolic manifolds. We say that a metric $g$ on $\cM^n$ is $\epsilon$-Einstein of order $\gamma$ if
\begin{equation}\label{ordered pinch}
\|h_g\|(x) \leq \epsilon \re^{- \gamma d(x, x_0)}.
\end{equation}
 
 \vskip 0.1in
 
\begin{theo}\label{maintheorem2}
Given $n\geq3$ and positive constants $k_0,k_1, v_0,\lambda_0$, and $\alpha$, let $\gamma > \frac 12 \alpha - \sqrt{\lambda_0}$. Then there exists $\epsilon>0$ depending only on $n, k_0, k_1, v_0, \lambda_0, \alpha, C_0$, and $\gamma$ such that the normalized Ricci
flow starting from a metric $g_0$ exists for all the time and converges exponentially to
an Einstein metric, provided that 
\begin{enumerate}[(1)]
\item
$\|Rm_{g_0}\|\leq k_0, \|\nabla Rm_{g_0}\|\leq k_1$
\item
the volume bound $\text{vol}(B_{g_0}(x, 1))  \geq v_0$, for all $x\in \cM^n$
\item
$g_0$ is with the non-degeneracy $\geq \lambda_0$
\item
$\int_{\cM^n} \exp( - \alpha d(x, x_0)) dv_{g_0} < C_0$  where $C_0$ is independent of $x_0$ 
\item
$g_0$ is $\epsilon$-Einstein of order $\gamma$.
\end{enumerate}
\end{theo}

We notice here that the volume condition (4) in Theorem \ref{maintheorem2} always holds for $\alpha > n-1$ on asymptotically hyperbolic manifolds $(\cM^n, \ g)$ (cf. Lemma \ref{lemm:analytic tool}). We then adopt a maximum principle (cf. Lemma \ref{maximum principle}) similar to  that in \cite{LT} \cite{EH} to show that the normalized Ricci flow, starting from an asymptotically hyperbolic metric, remains to be asymptotically hyperbolic and the conformal infinity at any finite time remains the same as the initial one. But, in order to maintain the decay rate of the Ricci pinching curvature $h$ near the space infinity of the limit metric at the time infinity of the normalized Ricci flow, the decay rate of the Ricci pinching curvature $h$ seems to be required to lie in the range 
$$
(\frac {n-1}2 - \sqrt{\frac{(n-1)^2}4 - 2}, \quad  \frac {n-1}2 + \sqrt{\frac{(n-1)^2}4 -2})
$$
for $n \geq 4$ (cf. \eqref{the decay range} in the proof of Theorem \ref{long time existence-asym}).  Therefore we have

\begin{theo}\label{maintheorem3}
Given $n \geq 5$ and positive constants $k_0,k_1, v_0,\lambda_0$, let $\gamma > 2$ and
$$
\gamma \in  (\frac {n-1}2 - \min\{\sqrt{\frac{(n-1)^2}4 - 2},  \sqrt{\lambda_0}\}, \quad  \frac {n-1}2 + \sqrt{\frac{(n-1)^2}4 -2}).
$$ 
Then there exists $\epsilon>0$ depending only on $n, k_0, k_1, v_0,\lambda_0$, and $\gamma$ such that the normalized Ricci
flow starting from an asymptotically hyperbolic metric $g_0$ of $C^2$ regularity remains to be asymptotically hyperbolic  of $C^2$ regularity for all the time and converges 
exponentially to a conformally compact Einstein metric of $C^2$ regularity with the same conformal infinity of $g_0$, provided that 
\begin{enumerate}[(1)]
\item
$\|Rm_{g_0}\|\leq k_0, \|\nabla Rm_{g_0}\|\leq k_1$
\item
the volume bound $\text{vol}(B_{g_0}(x, 1))  \geq v_0$, for all $x\in \cM^n$
\item
$g_0$ is with the non-degeneracy $\geq \lambda_0$
\item
$g_0$ satisfies 
$$ 
\|h_{g_0}\| (x) \leq \epsilon\re^{-\gamma d(x, x_0)} \  \rm{and} \
\|\nabla h_{g_0}\|\leq C\re^{-\gamma d(x, x_0)}.
$$
\end{enumerate}
\end{theo}

We like to point out that one may produce conformally compact Einstein metrics with less regularity if $\gamma \leq 2$. Notice that it is good enough to produce conformally compact Einstein metrics of $C^2$ regularity in the light of the regularity results in \cite{CDLS}. Even though it is not easy for now to produce conformally compact Einstein metrics using the normalized Ricci flow on asymptotically hyperbolic manifolds in general. It is nice to see that Theorem \ref{maintheorem3} can be used to prove perturbation existence results in \cite{GL} \cite{Lee} \cite{Bi}. 

\begin{theo} \label{maintheorem4}
Let $(\cM^n, \ g)$, $n\geq 5$,   be a conformally compact Einstein manifold of regularity $C^2$ with a smooth conformal infinity $(\partial\cM, \ [\hat g])$.  And suppose that the non-degeneracy of $g$ satisfies
\begin{equation}\label{nondegeneracy}
\sqrt \lambda>\frac {n-1}2 -2.
\end{equation}
Then, for any smooth metric $\hat h$ on $\partial\cM$, which is sufficiently
$C^{2,\alpha}$ close to  some $\hat g\in [\hat g]$ for any $\alpha\in (0,1)$, there
is a conformally compact Einstein metric on $\cM$ which is of $C^2$
regularity and with the conformal infinity  $[\hat h]$.
\end{theo}

The non-degeneracy condition \eqref{nondegeneracy} is assumed because our estimate of the decay of the curvature $h$ in time relies on the non-degeneracy and
the decay of the curvature $h$ in space at the initial time in this paper. But 
the non-degeneracy condition \eqref{nondegeneracy} is trivially satisfied for all non-degenerate conformally compact Einstein metric $g$ in dimension 5. More interestingly, Theorem \ref{maintheorem4} fully recovers the perturbation existence result in \cite{GL} for dimensions other than $4$ since the non-degeneracy of the hyperbolic n-space is $\frac {(n-1)^2}4$. Furthermore, when the non-degeneracy is small in dimensions higher than $5$, we can get the perturbation existences in \cite{Lee} \cite{Bi} by constructing the initial metric with higher decay rate of the Ricci pinching curvature for a conformal infinity that is close to the given one in $C^{k, \alpha}$ for appropriately large $k$ (please see Theorem \ref{general perturbation}).

Our paper is organized as follows: In Section 2 we use the ideas from \cite{Ye} \cite{LY} to obtain the decay of Ricci pinching curvature in time for finite time when the metrics along the normalized Ricci flow are non-degenerate and sufficiently Ricci pinched. In Section 3 we continue to use ideas from \cite{Ye} \cite{LY} to get global existences and convergences of normalized Ricci flows based on finite time results in the previous section. In Section 4 we introduce asymptotically hyperbolic manifolds and establish Theorem \ref{long time existence-asym}. The main tool is the maximum principle (Lemma \ref{maximum principle}) adopted from \cite{LT} \cite{EH}. The main challenge is to keep the regularity of conformal compactness at the time infinity. In Section 5 we recall the metric expansions in \cite{FG} for conformally compact Einstein metrics and apply normalized Ricci flows to reproduce perturbation existence results in \cite{GL} \cite{Lee} \cite{Bi}. We also calculate curvatures for asymptotically hyperbolic metrics and evolution equations for curvatures along the normalized Ricci flow in the appendices for the convenience of readers. 


\section{NRF of non-degenerate and Ricci-Pinched metrics}

For simplicity we will use NRF in short for the normalized Ricci flow from now on. Suppose that $(\mathcal{M}^n, \  g)$ is a complete and non-compact Riemannian manifolds. Let  $\Delta_L$ be the Lichnerowicz Laplacian with respect to the metric $g$ on symmetric
$2$-tensor $u$ as follows (c.f. \cite{Besse}):
\begin{equation}\label{Lichnerowiczlaplacian}
\Delta_L u_{ij}=-\Delta u_{ij} - 2R_{ipjq}u^{pq}+R_{iq}u_j^q + R_{jq}u_i^q,
\end{equation}
where $R_{ipjq}$, and $R_{ij}=g^{pq}R_{ipjq}$ are the
components of Riemann curvature tensor and Ricci curvature tensor of the metric $g$ respectively.
Note that in our notations, for any smooth function
$v$ on $(\mathcal{M},g)$, we have
$$
\Delta_L (v g_{ij})=-(\Delta v)g_{ij},
$$
where $\Delta $ is the Laplacian operator on function with
respect to metric $g$.  For convenience we make the following convention.

\begin{defi}\label{def.ofconditionB}
A complete  Riemannian manifold
$(\mathcal{M}^n, \  g)$ is said to satisfy the condition $B(k_0,  v_0, \lambda_0)$
with some positive constants $k_0,  v_0, \lambda_0$, if there
hold the following three conditions:
\begin{enumerate}[(i)]
\item bounded  curvature condition
$$
\|Rm(g)\|_g \leq k_0,
$$
\item the volume bound 
$$\text{vol}(B_{g_0}(x, 1))  \geq v_0$$
for all $x\in \cM^n$.
\item non-degenerate condition
$$
\int_ \mathcal{M}  \langle\big(\Delta_L+2(n-1)\big)u_{ij}, u_{ij}\rangle \geq \lambda_0 \int_
\mathcal{M}\|u\|^2
$$ 
holds for any symmetric $2$-tensor $u$ such that $\int_{\cM}(|\nabla u|^2+ |u|^2)dv < \infty$.
\end{enumerate}
\end{defi}

Due to the equivalence between Ricci flow and the normalized Ricci flow, the short time existence and some basic curvature estimates for the normalized Ricci flows have been established in Theorem 1.1 in \cite{Shi}. But for our purpose we would like to have the following curvature estimates, which can be proven via the maximum principle as in \cite{Shi}. Namely,

\begin{lemm}\label{curvature estimate}
Let $g_0$ be a Riemannian metric on $\cM^n$ which satisfies
$$\|Rm\|_{C^0(\cM)}\leq k_0, \|\nabla
Rm\|_{C^0(\cM)}\leq k_1$$ with some positive constants
$k_0, k_1$. Let $g(t)$ for $t\in [0, T]$ is the solution to \eqref{NRF} obtained in Theorem 1.1 in \cite{Shi}.  Then, for each $t\in(0, T]$, the following estimates hold on $\cM^n$:
\begin{eqnarray}\label{derivative estimate of curvature}
 &\|\nabla Rm\|_{C^0(\cM)}(t)\leq
C_2,\\\label{high derivative estimate of curvature} &\|\nabla^2
Rm\|_{C^0(\cM)}(t)\leq \frac{C_3}{\sqrt t},
\end{eqnarray}
where $C_2,C_3$ are some constants depending only on $n, k_0, k_1, T$.
\end{lemm}

\begin{proof}
Estimates \eqref{derivative estimate of curvature}  and \eqref{high derivative estimate of curvature} are included in (4) in Theorem 1.1 in \cite{Shi}. But it takes a more or less the same proof to verify them. For the convenience of the readers we briefly sketch the proof.  First we know from (4) in Theorem 1.1 in \cite{Shi} that
$$
\|Rm\|^2 (\cdot, t) \leq C_1
$$
for $t\in [0, T]$. To prove  \eqref{derivative estimate of curvature}
set
$$
\varphi_1=(a+\|Rm\|^2)\|\nabla
Rm\|^2-(a+k_0^2)\cdot k_1^2,
$$
where $a>0$ is some constant to be determined later.
It follows from the evolution equations (\ref{evolving equation 4}) and (\ref{evolving equation
5}) that $\varphi_1$ satisfies
$$
\fr{\pr}{\pr t}\varphi_1 \leq
\Delta\varphi_1- \fr{1}{4a^2}(\varphi_1+(a+k_0^2)k_1^2)^2
+C (\varphi_1+ (a+k_0^2)k_1^2).
$$
upon choosing $a= 7C_1$.

Arguing as in [Shi], for each point $x_0\in\cM$, we can find a
function $\xi$ such that $\xi(x)\equiv1$ in $B(x_0,1)$,
$\xi(x)\equiv0$ in $\cM\setminus B(x_0,2)$ and $0\leq\xi(x)\leq1$
for $x\in\cM$. Moreover, its derivatives satisfy the following
estimates:
$$
|\tilde{\nabla}\xi|^2_0 \leq  C\xi(x)\quad\quad\forall
x\in\cM
$$
$$
\tilde{\nabla}_\alpha\tilde{\nabla}_\beta\xi(x) \geq -C g_{\alpha\beta}(x)\quad\quad\forall
x\in\cM
$$
where $\tilde{\nabla}$ and $|\cdot|_0$ are those with respect to the metric $g(0)$.
Set the function $F_1(x,t)=\xi(x)\varphi_1(x,t)$,
$(x,t)\in\cM\times[0, T]$, and note that our assumptions imply
 $$F_1(x,0)\leq 0 \quad\quad\quad\quad\quad\quad \forall
x\in\cM.$$ If $F_1(x,t)\leq0$ for each $(x,t)\in\cM\times(0, T]$,
then the proof of (\ref{derivative estimate of
curvature}) is complete.  Otherwise there exists a point $(x_1,t_1)\in
B(x_0,2)\times(0, T]$ such that
$$F_1(x_1,t_1)=\max_{\cM\times[0, T]}F_1(x,t)>0.$$
Then
$$
\nabla F_1(x_1,t_1) = 0,
\fr{\pr F_1}{\pr t}(x_1,t_1) \geq 0,  \text{ and }   \Delta F_1(x_1,t_1) \leq 0,$$
which implies, at $(x_1, t_1)$,
\begin{equation}\label{Eq5}
\frac 1{4a^2}\xi\varphi_1^2\leq-\varphi_1\Delta\xi+C\varphi_1 + C,
\end{equation}
where one uses the fact that
$$
\nabla\varphi_1\cdot\nabla\xi \geq - C\varphi_1
$$
due to the properties of the function $\xi$ and $\nabla F=0$. Recall that
$$
-\Delta\xi=-g^{\alpha\beta}\tilde{\nabla}_\alpha\tilde{\nabla}_\beta\xi+g^{\alpha\beta}
(\Gamma^\gamma_{\alpha\beta}-\tilde{\Gamma}^\gamma_{\alpha\beta})\tilde{\nabla}_\gamma\xi,
$$
\begin{equation}\label{christoffel}
\frac {\partial\Gamma^\gamma_{\alpha\beta}}{\partial t} = - g^{\gamma\delta}(R_{\alpha\delta, \beta} + R_{\beta\delta, \alpha} - R_{\alpha\beta, \delta}),
\end{equation}
and (4) in Theorem $1.1$ of \cite{Shi}, for $t\in[0, T]$
$$
\|\nabla Rm \|\leq \fr{C}{\sqrt t}.
$$
We easily get
$$|\Gamma^\gamma_{\alpha\beta}(x_1,t_1)-\tilde{\Gamma}^\gamma_{\alpha\beta}(x_1)| \leq C \text{ and }
-\Delta\xi\leq C.
$$
Therefore we conclude from \eqref{Eq5} that
$$
\frac 1{4a^2}F_1^2\leq C F_1+ C
$$
which implies that
$$
\|\nabla Rm\| (\cdot, t) \leq C(n,k_0,k_1, T).
$$
for all $t\in [0, T]$.

To prove \eqref{high derivative estimate of curvature}, we  set
$\varphi_2=(b+ \|Rm\|^2 + \|\nabla Rm\|^2)\|\nabla^2 Rm\|^2$ and
choose suitable $b$ such that
$$ \fr{\pr}{\pr t}\varphi_2\leq
\Delta\varphi_2-C \varphi_2^2+C \varphi_2 + C.$$
Considering $F_2= t\xi\varphi_2$, one similarly obtains the estimate \eqref{high derivative estimate of curvature}
\end{proof}

An immediate consequence of \eqref{high derivative estimate of curvature} is that the sectional curvature point-wisely changes
little in a short time period along NRF. Next  we specify the Ricci curvature pinching conditions in this paper.

\begin{defi}\label{Einstein}
A Riemannian metric $g$ on $\cM^n$ is called $\varepsilon$-Einstein
if it satisfies on $\cM$ that
$$\sup_{x\in \cM^n}\|Ric(g)+(n-1)g\|_g(x) \leq \varepsilon.$$
A metric $g$ is said to be $\varepsilon$-Einstein of
order $\gamma$ if, for a positive number $\delta$, there holds for any $x\in\cM$ that
$$\|Ric(g)+(n-1)g\|_g(x)\leq \varepsilon e^{-\gamma d(x,x_0)},$$ where $d(x,x_0)$ is the
distance from $x$ to some fixed point $x_0$ with respect to $g$.
\end{defi}

Let us first investigate the behavior of the tensor $h$ along
NRF when the metrics along the flow are known to be
$\varepsilon$-Einstein and satisfy the condition $B(k_0,  v_0,
\lambda_0)$. This is easy, particularly after \cite{Ye} \cite{LY}, yet the important initial step. We adopt the approach in \cite{Ye} to obtain $C^0$
decay estimates of $h$ from $L^2$ decay estimates via Moser iteration. To derive
$L^2$ decay estimate from the non-degenerate condition on complete and non-compact manifold,
we follow the interesting heat kernel estimates of
Grigor$'$yan \cite{Gri} given in \cite{LY}.

\begin{lemm}\label{lemma3}
Let $g(t), t\in[0,T]$, be a solution to \eqref{NRF}. We assume that $g(t)$ satisfies
$$
\|g(\cdot, t) - g(\cdot, 0)\|_g \leq \varepsilon,
$$
the condition $B(k_0, v_0,\lambda_0)$, and $\varepsilon$-Einstein for each $t\in [0, T]$ and sufficiently small $\varepsilon$. And we assume that
$$
\int_{\cM^n}\exp\big(-\alpha d(x,x_0)\big)\mathrm{dv_{g(0)}}<\infty
$$
for a positive number $\alpha$. In addition, suppose that
$$
\int_{\cM} \|h\|^2  \mathrm{dv_{g_0}}<\infty.
$$
Then for any $(x,t)\in\cM\times[\tau,T]$, we have
\begin{enumerate}[(i)]
\item
\begin{equation*}
\|h\|^2
\leq C\re^{-(2\lambda_0 -C \varepsilon)t}
\int_{\cM} \|h\|^2  \mathrm{dv_{g_0}}
\end{equation*}
\item
 \begin{equation*}
\|\nabla h \|(x,t)+\|\nabla^2 h\|(x,t) \leq
C\re^{-(2\lambda_0 -C \varepsilon)t}
\int_{\cM} \|h\|^2  \mathrm{dv_{g_0}}
\end{equation*}
\end{enumerate}
where $C$ is a constant depending only on $\tau,k_0, \delta_0, \lambda_0, n$.
\end{lemm}

\begin{proof} The first step is the same as observed in \cite{Ye}.
Since the sectional curvature with $g(t)$ is
assumed to be bounded by $k_0$ the evolution equation (\ref{evolving equation
2}) yields
$$
\frac{\partial}{\partial t}\|h\|^2\leq \Delta \|h\|^2+ C \|h\|^2.
$$
We also know that the Sobolev constant is bounded because the curvature is bounded and the injectivity radius is bounded from below from the condition $B(k_0, v_0, \lambda_0)$ (cf. Lemma \ref{injectivity radius} and \cite{Che} \cite{Heb}). Hence the standard Moser iteration in the parabolic ball
$B_0(x,\sqrt{\tau})\times[t-\tau,t]$ implies
\begin{equation}\label{eq23}
\sup_{B_0(x,\sqrt{\frac{\tau}{2}})\times[t-\frac{\tau}{2},t]}||h||^2
\leq
C(n,\tau,k_0)\int_{t-\tau}^{t}\int_{B_0(x,\sqrt{\tau})}
\|h\|^2(y,s)\mathrm{d}y\mathrm{d}s.
\end{equation}

In the second step we follow the idea in \cite{LY} to derive the decay from the non-degeneracy of $g(t)$. To do so we recall the auxiliary function from \cite{LY}
$$
\xi(y,s)=-\fr{d_0^2(y)}{(2+C_0\varepsilon)(t-s)},
$$
where $d_0(y)$ is the distance function from the point $y$ to the
geodesic ball $B_0(x,\sqrt\tau)$ with respect to the initial metric $g_0$ and
$C_0$ is chosen so that \be\label{propertyofxi}
\xi_s+\fr{1}{2}\|\nabla\xi\|^2\leq 0. \ee 
We then set
$$
J(s)=\int_{\cM^n}\|h\|^2(y,s)\re^{\xi(y,s)}\mathrm{d}y.
$$
With the volume growth condition and the fact that $g(t)$ are all  quasi-isometric to $g_0$ one sees that $J(s)$ is finite for all $0<s<t<T$.
The important observation here is that the non-degeneracy implies the exponential decay of $J(s)$ when it evolves.  We compute
\begin{eqnarray*}
\frac{\ud{J}}{\ud s}(s)&=& \int_{\cM}2\langle\Delta h_{ij}-2R_{ipjq}h_{pq}-2h_{ip}h_{pj},
h_{ij}\rangle\re^{\xi}\\
&&+\|h\|^2\re^{\xi}\xi_s\ud y+C \varepsilon J(s),\nonumber
\end{eqnarray*}
where we used the evolution equation (\ref{evolving equation 1}) and the
$\varepsilon$-Einstein condition for $g(t)$.
It is crucial to realize from (\ref{propertyofxi}) that
\begin{equation}\label{Eq12}
\int_{\cM}2(\Delta h)h\re^{\xi}+\|h\|^2\re^{\xi}\xi_s\ud
y \leq \int_{\cM}2<\Delta (\re^{\fr{\xi}{2}}h), \ (\re^{\fr{\xi}{2}}h)>
\end{equation}
and therefore
\begin{equation}\label{Eq14}
\frac{\ud{J}}{\ud s}(s) \leq
-2\int_{\cM} <(\Delta_L +2(n-1))(\re^{\frac \xi 2}h), (\re^{\frac \xi 2}h)> + C\varepsilon J(s),
\end{equation}
which implies
\begin{equation}\label{nondegenerate}
\frac{\ud{J}}{\ud s}(s) \leq
-(2\lambda_0 - C \varepsilon)J(s)
\end{equation}
by the non-degeneracy property of $g(t)$ and the fact that
$$
\int_{\cM} (\|\nabla (\re^{\frac \xi 2}h)\|^2 + \re^\xi \|h \|^2) < \infty
$$
in the light of curvature bounds and the volume growth assumption. Therefore we obtain
\begin{equation}\label{decay of J}
J(s)\leq \re^{-(2\lambda_0 - C\varepsilon )s} J(0).
\end{equation}
Going back to \eqref{eq23} we thus derive that
\begin{equation}\label{decay-h}
\sup_{B_0(x,\sqrt{\frac{\tau}{2}})\times[t-\frac{\tau}{2},t]}
\|h\|^2 \leq C \re^{- (2\lambda_0 - C \varepsilon) t}\int_{\cM} \|h\|^2 \re^{\xi(\cdot, 0)} \mathrm{dv_{g(0)}}
\end{equation}
which implies (i) of the conclusion. And the conclusion (ii) follows from (i) and some standard estimates.
\end{proof}

When the initial metric is $\varepsilon$-Einstein of order $\gamma$, instead of having $L^2$ of the curvature $h$ finite, with a little changes of the proof, we conclude that

\begin{lemm}\label{lemma:pinch with order} Let $g(t), t\in[0,T]$, be a solution to \eqref{NRF}. We assume that $g(t)$ satisfies 
$$
\|g(\cdot, t) - g(\cdot, 0)\|_g \leq \varepsilon,
$$
the condition $B(k_0, v_0,\lambda_0)$, and $\varepsilon$-Einstein, for each $t\in [0, T]$ and sufficiently small $\varepsilon$.
And we assume that
\begin{equation}\label{to be proved}
\int_{\cM} \exp(-\alpha d(x, x_0))dv_g\leq C
\end{equation}
for any $x_0\in\cM$, where $C$ is independent of $x_0$. In addition, suppose that $g_0$ is $\epsilon_0$-Einstein of order $\gamma$ such that
\begin{equation}\label{delta}
\gamma > \frac 12(\alpha - 2\sqrt{\lambda_0}).
\end{equation}
Then  for any $(x,t)\in\cM\times[\tau,T]$, we have
\begin{enumerate}[(i)]
\item
\begin{equation*}
\|h\|^2
\leq C\epsilon_0 \re^{- \lambda_1t}
\end{equation*}
\item
 \begin{equation*}
\|\nabla h \|(x,t)+\|\nabla^2 h\|(x,t) \leq
 C\epsilon_0 \re^{- \lambda_1 t}\end{equation*}
\end{enumerate}
for any $0 < \lambda_1 < 2\lambda_0$.
\end{lemm}

\begin{proof} This is simply because,  in (\ref{decay-h}) in the previous proof, one notices that
$$
\frac {d_0^2}{((2 +C\varepsilon)t} + (2\lambda_0 - C\varepsilon - \lambda_1)t \geq 2 \sqrt{\frac {2\lambda_0 - C\varepsilon - \lambda_1}{2 + C\varepsilon}} d_0
$$
and
$$
\int_{\cM} \exp(-ad(y, x))\exp(-bd(y, x_0)) dv_g \leq 2C
$$
when $a+b \geq\alpha$ in the light of the volume assumption \eqref{to be proved}.
\end{proof}

We remark that the volume condition \eqref{to be proved} will be proven to hold on asymptotically hyperbolic manifolds in Lemma \ref{lemm:analytic tool}.
In particular, in the cases of perturbation of hyperbolic space, as discussed in \cite{LY}, then $\alpha > n-1$ and $\lambda \geq
\frac {(n-1)^2}{4}-C\varepsilon$. Therefore \eqref{delta} becomes simply $\gamma > 0$ in \cite{LY}.


\section{Long time existence and convergences}

In this section we show that NRF starting from a non-degenerate metric with a sufficiently pinched Ricci curvature will remain to be non-degenerate and sufficiently pinched in Ricci curvature. Therefore we will be able to use Lemma \ref{lemma3} and Lemma  \ref{lemma:pinch with order} to prove long time existence and convergence theorems on complete non-compact manifolds. This approach is more or less the same as the one taken in \cite{Ye}. The key is to investigate how curvature bound, the volume $\text{vol}(B(x, 1))$ lower bound  and non-degeneracy evolve along NRF.  The strategy is to first use curvature bounds of Theorem 1.1 in \cite{Shi} and Lemma \ref{curvature estimate} in the previous section to make sure those quantities change little in a very short time, then using exponential decay estimates of Lemma \ref{lemma3} and Lemma  \ref{lemma:pinch with order} to keep those quantities change little even after an arbitrarily  long time.

It is known that one can bound the injectivity radius from below provided
that one has volume lower bound and bounded curvature (cf. \cite{Che}).  For the convenience of readers we recall that

\begin{lemm}\label{injectivity radius}
Let $(\cM^n, \ g)$ be a complete Riemannian manifold and that
$$
|Sec_{\cM}|\leq k \quad\text{and}\quad volB(x,1)\geq v
$$ 
for all $x\in\cM^n$, where $B(p, 1)$ is the unit geodesic ball at $p$. Then 
\begin{equation}\label{injectivity radius}
inj_{\cM} \geq \delta (n, k, v).
\end{equation}
\end{lemm}

Let us solve the short time problem first.

\begin{lemm}\label{short time}
Assume that $g_0$ satisfies the condition
$B(k_0,v_0,\lambda_0)$ and is $\varepsilon_0$-Einstein.  In addition we assume that $|\nabla Rm|\leq k_1$ for $g_0$. Then, for any $\epsilon>0$, 
there exists a positive number $\tau$ depending only on
$k_0,k_1,\lambda_0, n$, and $\epsilon$ such that the solution $g(t)$ to \eqref{NRF} satisfies the conditions
$B(k_0+\epsilon, v_0 - \epsilon, \lambda_0-\epsilon)$ and $(\varepsilon_0+\epsilon)$-Einstein for all $t\in[0,\tau]$.
\end{lemm}

\begin{proof} The proof is to verify the continuity of all the quantities involved in the conditions $B(k_0, v_0, \lambda_0)$ and $\varepsilon$-Einstein.  In the light of Lemma \ref{curvature estimate} we easily see that there is $\tau >0$ such that
$$
\|Rm(g(t))\|_{g(t)} \leq k_0+\epsilon \text{ and } \|h(g(t))\|_{g(t)} \leq \varepsilon_0+\epsilon
$$
for all $t\in [0, \tau]$. 

For the non-degeneracy, we recall the quadratic form $Q$ on symmetric $2$-tensors:
\begin{equation}\label{q-form}
\aligned
Q(u,g) & = \int_{\cM} <(\Delta_L +2(n-1))u, \ u> dv_g \\
& = \int_{\cM^n}\|\nabla u\|^2+2\int_{\cM^n}R_{ipjq}u^{pq}u^{ij}
+\int_{\cM^n}(h_{ik}u_{j}^{k}+h_{jk}u_{i}^{k})u^{ij}
\endaligned
\end{equation}
and the definition of  non-degeneracy:
$$
\lambda (g)=\inf\frac{Q(u,g)}{\int_{\mathcal{M}} \|u\|^2 _g dv_g},
$$
where  infimum is taken over all symmetric $2$-tensor $u$ on $(\mathcal{M},g)$ with
$$
\|u\|^2 _{W^{1,2}(\mathcal{M},g)}=\int_{\mathcal{M}}(\|\nabla u\|^2
+\|u\|^2)dv_g <\infty,
$$
and $\nabla$ and norm $\|\cdot\|$ are taken with respect to metric $g$.  It is then straightforward to estimate
$$
\lambda (g(t)) \geq \lambda_0-\epsilon
$$
when $\|g(t)  -g_0\|,  \|\Gamma^\gamma_{\alpha\beta}(g(t)) - \Gamma^\gamma_{\alpha\beta}(g_0)\|$, and $ \|Rm(g(t))- Rm(g_0)\|$ are all sufficiently small for all $t\in [0, \tau]$.

Finally it is easily seen that
$$
\text{vol}(B_{g(t)}(x, 1)) \geq v_0 - \epsilon
$$ 
when $\tau$ is sufficently small.
\end{proof}

Next let us show that NRF starting from a non-degenerate and Ricci sufficiently pinched metric will remain non-degenerate and Ricci pinched forever. The idea of the proof is the same as the one in \cite{Ye}.

\begin{lemm}\label{long time}
Suppose $(\mathcal{M}, \ g_0)$ is a complete Riemannian manifold
satisfying the conditions $B(k_0, v_0, \lambda_0)$, $|\nabla Rm|\leq k_1$,  and
$\varepsilon_0$-Einstein, where $8C\varepsilon_0 < \lambda$ in Lemma \ref{lemma3} for the metric $g(t)$ satisfying the
condition $B(2k_0, \frac 12 v_0, \frac 12\lambda_0)$ along
NRF. And suppose that
$$
\int_{\cM} \exp( - \alpha d(x, x_0)) dv_{g_0} < \infty
$$
for some positive number $\alpha$. Then, for any $\epsilon >0$,  there is a small number
$\delta > 0$ such that, the solution to NRF exists in $[0,
+\infty)$ and satisfies the conditions $B(k_0 + \epsilon, v_0 - \epsilon,
\lambda_0 -\epsilon)$ and $(\varepsilon_0+\epsilon)$-Einstein for all $t\in
[0,+\infty)$, provided that
$$
\int_{\cM} \|h\|^2 dv_{g_0} \leq \delta.
$$
\end{lemm}

\begin{proof} First, for NRF starting from a metric $g_0$ satisfying the conditions $B(k_0, v_0, \lambda_0)$ and $|\nabla Rm|\leq k_1$, we apply Lemma \ref{short time}
and find a number $\tau >0$ such that $g(t)$ satisfies the condition
$B(k_0+ \frac 12\epsilon, v_0 - \frac 12\epsilon, \lambda_0-\frac 12\epsilon)$
and is $(\varepsilon_0+\frac 12\epsilon)$-Einstein for each $t\in [0, \tau]$.
It is important that $\varepsilon_0$ is small enough  that
$8C\varepsilon_0 < \lambda_0.$ for the metrics $g(t)$ in $B(2k_0, \frac 12v_0, \frac 12\lambda_0)$ in Lemma \ref{lemma3}.

Now, for any $T > \tau$, if $g(t)$ satisfies the condition $B(k_0+ \frac 34\epsilon, v_0- \frac 34\epsilon, \lambda_0-\frac 34\epsilon)$ and is $(\varepsilon_0+\frac 34\epsilon)$-Einstein for all $t\in [0, T]$, based on the estimates (i) and (ii) in Lemma \ref{lemma3} and the evolution equations \eqref{NRF} \eqref{evolvingequation 3}  \eqref{christoffel} we have
$$
\|g(t) - g(\tau)\| \leq C (\int_{\cM}\|h\|^2 dv_g)^{\fr 1 2}
$$
$$
\|\Gamma^\gamma_{\alpha\beta}(g(t)) -
\Gamma^\gamma_{\alpha\beta}(g(\tau))\| \leq  C (\int_{\cM}\|h\|^2
dv_g)^{\fr 1 2}
$$
$$
\|Rm(g(t)) - Rm(g(\tau))\| \leq C (\int_{\cM}\|h\|^2dv_g)^{\fr 1 2}
$$
where $C$ depends only $\tau, k_0, v_0, \lambda_0, n$ for $t\in
[\tau, T]$. Therefore the Riemann curvature bound, injectivity
radius, non-degeneracy, and the Ricci pinching all change as small
as desired, independently of $T$,  provided $\int_{\cM}\|h\|^2dv_g$
is sufficiently small. For instance we take $\int_{\cM}\|h\|^2 dv_g<
\epsilon$ so that $g(t)$ remains to satisfy the conditions $B(
k_0+\frac 58\epsilon, v_0-\frac58\epsilon, \lambda_0-\frac 58\epsilon)$ and
$(\varepsilon_0+\frac58\epsilon)$-Einstein for all $t\in [0, T]$. What is
important here is that the choice of $\epsilon$ is independent of
$T$ because of the exponential decay in Lemma \ref{lemma3}.

Thus we may argue as follows: Suppose otherwise we have a largest number $T_0 < \infty$ such that $g(t)$ exists and satisfies the condition $B(k_0+\frac 34\epsilon, 
v_0-\frac34\epsilon, \lambda_0-\frac34\epsilon)$ and $(\varepsilon_0+\frac34\epsilon)$-Einstein for all $t\in [0, T_0)$. Then, the above argument, which is based on 
Lemma \ref{lemma3}, tells us that in fact
$g(t)$ satisfies the conditions 
$B(k_0+\frac 58\epsilon, v_0-\frac58\epsilon, \lambda_0-\frac 58\epsilon)$ and
$(\varepsilon_0+\frac58\epsilon)$-Einstein for all $t\in [0, T_0]$, which contradicts with the maximality of $T_0$.
 \end{proof}

Using a similar argument, based on Lemma \ref{lemma:pinch with order}, we therefore also have

\begin{lemm}\label{long time-asym}
Suppose that $(\mathcal{M}, \ g_0)$ is a complete Riemannian manifold
satisfying the conditions $B(k_0, v_0, \lambda_0)$ and $|\nabla Rm| \leq k_1$,  and that
$$
\int_{\cM} \exp( - \alpha d(x, x_0)) dv_g \leq C
$$
for some positive number $\alpha$ and $C$, where $C$ is independent
of $x_0$. Then, for any number $\epsilon >0$ and a number $\gamma > \frac 12(\alpha - 2\sqrt \lambda_0)$, there is a small number
$\delta > 0$ such that, the solution to NRF exists in $[0,
+\infty)$ and satisfies the conditions $B(k_0+\epsilon, v_0-\epsilon,
\lambda_0-\epsilon)$ and $(\varepsilon_0+\epsilon)$-Einstein for all $t\in
[0, +\infty)$, provided that the initial metric $g_0$ is
$\delta$-Einstein of the order $\gamma$.
\end{lemm}

Finally we state and prove main theorems in this section including the convergences for the above global NRF at the infinity of the time.

\begin{theo}\label{long time existence}
For any $n\geq3$ and positive constants $k_0, k_1, v_0, \lambda_0, \alpha$, there exists a positive number $\delta>0$ depending only on $n, k_0, k_1, v_0, \lambda_0$ such that NRF on $\cM^n$ starting from a metric $g_0$ satisfying the conditions $B(k_0,v_0,\lambda_0)$, $|\nabla Rm| \leq k_1$, and
$$
\int_{\cM} \exp (- \alpha d(x, x_0)) dv_g < \infty
$$
exists for all the time  if the initial metric is $\delta$-Einstein and
\be\label{L2condition}
\int_{\cM}\|h\|_g^2 dv_g \leq \delta.
\ee
Moreover NRF converges exponentially to a
non-degenerate Einstein metric $g_\infty$ in $C^\infty$ as $t\to\infty$.
\end{theo}

\begin{proof} First, from Lemma \ref{long time}, we pick up positive numbers $\varepsilon_0$ and $\delta < \varepsilon_0$ such that NRF from a metric $g_0$ satisfying
\begin{enumerate}[(1)]
\item
the conditions $B(k_0, v_0, \lambda_0)$ and $|\nabla Rm|\leq k_1$;
\item
$\varepsilon_0$-Einstein;
\item
 \begin{equation*}
\int_{\cM} \|h\|^2_g  dv_g \leq \delta,
\end{equation*}
\end{enumerate}
exists and satisfies the condition $B(2k_0, \frac 12v_0,
\frac 12\lambda_0)$ and $2\varepsilon_0$-Einstein for all the
time. Therefore
$$
\|h\|^2(\cdot, t) \leq C \re^{ -(\lambda_0-C \epsilon) t}
$$
for all the time big. Thus NRF converges exponentially to a non-degenerate Einstein metric $g_\infty$ as $t\to\infty$. The NRF actually converges to $g_\infty$ in $C^\infty$ due to the curvature estimates and the lower bound on the injectivity radius.
\end{proof}

Similarly we also can conclude

\begin{theo}\label{long time existence-pre-asym}
Given any $n\geq3$ and positive constants $k_0, k_1, v_0, \lambda_0, \alpha$, let $\gamma > \frac \alpha 2 - \sqrt \lambda_0$. Suppose that a given metric $g_0$ on $\cM^n$ satisfies the conditions $B(k_0, v_0,\lambda_0)$, $|\nabla Rm|\leq k_1$,  and
$$
\int_{\cM} \exp( -\alpha d(x,  x_0))dv_g \leq C_0,
$$
where $C$ is independent of $x_0$. Then there exists a positive
number $\delta >0$ depending only on
$n, k_0, k_1, v_0, \lambda_0, \alpha, \gamma, C_0$ such that NRF on $\cM^n$ starting from the metric $g_0$ exists for all the time, provided that $g_0$ is $\delta$-Einstein of order $\gamma$. Moreover NRF converges exponentially to a
non-degenerate Einstein metric $g_\infty$ in $C^\infty$ as $t\to\infty$.
\end{theo}


\section {NRF on asymptotically hyperbolic manifolds}

In this section we consider NRF on asymptotically hyperbolic manifolds. To see how NRF starting from an asymptotically hyperbolic metric remains to be asymptotically hyperbolic with the same regularity, at least up to certain order, we need to estimate the decay rate of $h$ at infinity of manifold $(\cM, g(t))$.  Our approach is to rely on a maximum principle.

Let us first introduce asymptotically hyperbolic manifolds. Suppose that $\cM^n$ is a smooth manifold with boundary $\partial\cM^{n-1}$. A defining function of the boundary is a smooth function $x: \bar\cM \to R^+$ such that, 1) $x > 0$ in $\cM$; 2) $x = 0$ on $\partial\cM$; 3) $dx\neq 0$ on $\partial\cM$. A metric $g$ on $\cM$ is said to be conformally compact if $x^2g$ is a Riemannian metric on $\bar\cM$ for a defining function $x$. The metric $g$ is said to be conformally compact of regularity $C^{k, \alpha}$ is $x^2g$ is a $C^{k, \alpha}$ metric on $\bar\cM$. The metric $\bar g = x^2 g$ induces a metric $\hat g$ on the boundary $\partial\cM$ and the metric $g$ induces a conformal class of metric $[\hat g]$ on the boundary $\partial\cM$ when defining functions vary. The conformal manifold $(\partial\cM, \ [\hat g])$ is called the conformal infinity of the conformally compact manifold $(\cM, \  g)$.

$(\cM, \ g)$ is said to be asymptotically hyperbolic if it is conformally compact and the sectional curvature of $g$ goes to $-1$ approaching the boundary at the infinity. The most basic and important fact about asymptotically hyperbolic manifolds is that, for a choice of a representative $\hat g\in [\hat g]$, there is a unique so-called geodesic defining function $r$ such that there is a coordinate neighborhood of the infinity $(0, r_0)\times\partial\cM \subset \cM$ where the metric $g$ is in the following normal form
\begin{equation}\label{normal form}
g = r^{-2} (dr^2 + g_r)
\end{equation}
where $g_r$ is one-parameter family of metrics on $\partial\cM$ and $(g_r)|_{r=0} = \hat g$.

Our goal is to consider Theorem \ref{long time existence-pre-asym} in the context of asymptotically hyperbolic manifolds. For that purpose we need a basic property of asymptotically hyperbolic manifolds regarding \eqref{to be proved}. A very efficient way to use the asymptotic hyperbolicity of asymptotically hyperbolic manifolds is the boundary M\"{o}bius charts introduced in \cite{Lee} (cf. Lemma 6.1 in \cite{Lee}).

\begin{lemm}\label{lemm:analytic tool} Suppose that $(\cM^n, \ g)$ is an asymptotically hyperbolic manifold. Then
$$
\int_{\cM} \exp(-\alpha d(x, x_0)) dv_g \leq C
$$
for any constant $\alpha > n-1$, where $C$ is independent of $x_0$.
\end{lemm}\begin{proof} First of all this lemma is obviously true on the hyperbolic space form simply because the isometry group of the hyperbolic space form is transitive.

For any given $x_0\in \cM$ it is clear that
$$
\int_{\cM} \exp(-\alpha d(x, x_0)) dv_g < \infty
$$
for any number $\alpha > n-1$. The issue is to prove that $C$ can be uniform for all $x_0\in \cM$. Hence we will only need to deal with the cases when $x_0$ is near the boundary at the infinity. In other words, we may assume $x_0$ is in boundary M\"{o}bius charts  $U_1 = \{ (r, \theta): r\in (0, r_1) \text{ and } |\theta| < r_1\} \subset U_2 = \{(r, \theta): r \in (0, 2r_1) \text{ and } |\theta| < 2r_1\}$. Using Lemma 6.1 in \cite{Lee} on the boundary M\"{o}bius charts it is clear that
$$
\int_{U_2} \exp(-\alpha d(x, x_0))dv_g \leq C
$$
for some constant independent of $x_0$. Therefore it suffices to show
$$
\int_{\cM\setminus U_2}\exp(-\alpha d(x, x_0))dv_g \leq C
$$
for some constant independent of $x_0$. In fact we may just let
$$
\int_{\{x\in \cM: r \geq r_1\}} \exp (-\alpha d(x, x_0))dv_g \leq \text{vol} (\{x\in \cM: r \geq r_1\}).
$$
Now for $x\in \{r \in (0, r_1)\}\setminus U_2$ we notice that $d(x, x_0) \geq  \frac {r_1}r$. Because the geodesic realizing the distance from $x$ to $x_0$ must enter $U_1$ and the distance from $x\in \{r \in (0, r_1)\}\setminus U_2$ to the boundary of $U_1$ is at least $\frac {r_1}r$. Therefore
$$
\int_{\{r \in (0, r_1)\}\setminus U_2} \exp(-\alpha d(x, x_0))dv_g \leq C \int_0^{r_1} \int_{\partial\cM} r^{ -n}\exp( -\alpha \frac {r_1}r) d\sigma dr\leq C.
$$
\end{proof}

In \cite{EH}, the authors adapted a method from \cite{LT} and extended the maximum principle on non-compact manifolds in \cite{LT}
to allow the metrics to be time dependent.  For our purpose we will need a variant of Theorem 4.3 in \cite{EH}. To clarify the terminology, a complete Riemannian manifold $(\cM, \ g)$ with boundary here is a non-compact manifold with a compact boundary $\partial\cM$. The typical example is the exterior of ball in the Euclidean space. Also $u_+ = \max\{u, 0\}$ as usual.

\begin{lemm}\label{maximum principle}
Suppose that $(\cM^n,g(t))$ is a smooth family of complete Riemannian manifolds with boundary $\partial \cM$ for $t\in [0, T]$.
Let $u$ be a function on
$\cM\times[0,T]$ which is smooth on  $\mathring{\cM}\times(0,T]$ and continuous
on $\cM\times[0,T]$. Assume that $u$ and $g(t)$ satisfy
\begin{enumerate}[(i)]
\item the differential inequality
$$\frac{\partial}{\partial
t}u-\Delta_{g_t}u\leq \textbf{a}\nabla u+ b u,$$ where the vector
$\textbf{a}$ and the function $b$ are uniformly bounded
$$\sup_{\cM\times[0,T]}|\textbf{a}|\leq\alpha_1,
\qquad \sup_{\cM\times[0,T]}|b|\leq\alpha_2,$$ with some constants
$\alpha_1,\alpha_2<\infty.$
\item $$\sup_{\cM}u(x,0)\leq 0 \text{ and } \sup_{\partial\cM \times[0, T]} u(x, t) \leq 0.$$
\item $$ \int_0^T\int_{\cM^n} \exp{[-\alpha_3 d^t(y, p)^2]} u_+^2 (y)dv_t(y)dt < \infty$$ for some positive number $\alpha_3$.
\item $$\sup_{\cM\times[0,T]}|\frac{\partial}{\partial
t}g(x,t)|\leq \alpha_4$$ with some constant $\alpha_4<\infty$.
\end{enumerate}
Then we have $u\leq 0$ on $\cM\times[0,T]$.
\end{lemm}

\begin{proof} The proof only differs from the proof in \cite {EH} in the last step. For the convenience of readers we sketch it and readers are referred to the proof of Theorem 4.3 in \cite{EH} for details. First we use the same function
$$
z (y, s) = - \frac {d_s^2(p, y)}{16(2\eta_0 - s)},
$$
where $d_s(p, y)$ is the distance function with respect to the
metric $g(s)$ and $\eta_0 < \min\{T, \frac 1{ 4\alpha_3}\}$.
Secondly we consider the same function
$$
u_K = \max\{\min\{u, K\},  0\}.
$$
Thirdly the cut-off function $\phi$ is defined as usual: $\phi = 1$ on $ B_p(R)$ ; $\phi = 0$ outside of $B_p(R+1)$, and $\phi$ is time independent, as in \cite{EH}. Following the same calculations in \cite{EH} we arrive at (please see the equation at the bottom of page 565 of \cite{EH})
\begin{equation}\label{eh-lt}
\aligned
e^{-\beta s}\int_{B_p(R)}  e^z u_K^2|_{s=\eta} &  \leq 4 \int_0^\eta e^{-\beta s} \int_{\cM_s\cap B_p(R+1)} e^z (|\nabla u|^2 - |\nabla u_K|^2)  \\
& \quad  + C \int_0^\eta e^{-\beta s}\int_{B_p(R+1)\setminus B_p(R)} e^z u_K^2 \endaligned
\end{equation}
for any $\eta \in (0, \eta_0]$ and $\cM_s = \{x\in \cM: u(x, s) > 0\}$. Notice that $\cM_s\cap\partial\cM = \emptyset$ for all $s\in [0, T]$ in the light of the assumption $(ii)$. Because of the choice of $\eta_0$ we know from the assumption $(iii)$
$$
\int_0^\eta\int_{\cM} e^z u_+^2 < \infty.
$$
Therefore, instead of asking $R\to\infty$ first to get rid of the
second term, here we ask $K\to\infty$ first to get rid of the first
term.  Since $|\nabla u_K|^2 \to |\nabla u|^2$ monotonically on
$\cM_s$ as $K\to\infty$ we have from \eqref{eh-lt}
$$
e^{-\beta s}\int_{B_p(R)} e^z (u_K)^2|_{s=\eta}  \leq  \quad   C \int_0^\eta e^{-\beta s}\int_{B_p(R+1)\setminus B_p(R)} e^z u_+^2.
$$
Thus, from the assumption (iii), taking $R\to\infty$ we finally reach
$$
\int_{\cM} e^z (u_+)^2 |_{s=\eta}\leq 0
$$
which implies that $u(\cdot, \eta)\leq 0$ for all $\eta\in [0, \eta_0]$. Then inductive argument implies that $u(\cdot, s)\leq 0$ for all $s\in [0, T]$.
\end{proof}

First let us state and prove an easy consequence of Lemma
\ref{maximum principle} and obtain the decay estimates that are
time-dependent, in other words, short time regularity for $g(t)$ as asymptotically hyperbolic metrics. 

\begin{lemm}\label{t-dependent} Suppose that $g(t)$, $t\in [0, T]$, is NRF starting from an asymptotically hyperbolic metric $g$ satisfying
$\|Rm\|\leq k_0$ and $\|\nabla Rm\|\leq k_1$. Then there exist numbers $C$, depending on $k_0, k_1, n, C_0$, and $T$ such that
$$
\|h\|(\cdot, t) \leq C r^\gamma, \ \|\nabla h\| (\cdot, t) \leq C r^\gamma \text{ and } \|\nabla^2 h\|(\cdot, t) \leq  \frac C{\sqrt t} r^\gamma
$$
for all $t\in [0, T]$, if
$$
\|h\|(\cdot, 0) \leq C_0r^\gamma \text{ and } \|\nabla h\|(\cdot, 0)\leq C_0 r^\gamma.
$$
\end{lemm}

\begin{proof} In the light of the curvature estimates in Lemma \ref{curvature estimate},  from
the evolution equations (\ref{evolving equation
2}),(\ref{evolvingequation 7}),(\ref{evolvingequation 8}), we obtain
\begin{eqnarray*}
\frac{\partial}{\partial t}\|h\|^2&\leq&\Delta \|h\|^2-2\|\nabla
h\|^2+C||h||^2,\\
 \frac{\partial}{\partial
  t}\|\nabla h\|^2&\leq&\Delta\|\nabla h\|^2-2\|\nabla^2 h\|^2
  + C(\|h\|^2+\|\nabla h\|^2),\\
  \frac{\partial}{\partial
  t}(t\|\nabla^2 h\|^2)&\leq&\Delta(t\|\nabla^2 h\|^2)-2t\|\nabla^3 h\|^2\\
  &&+(1+Ct)\|\nabla^2 h\|^2+C(||h||^2+||\nabla h||^2),
\end{eqnarray*}

Let $r$ be a fixed geodesic defining function of the asymptotically
hyperbolic metric $g_0$, one knows the fact that $|\Delta_g r|\leq Cr$
and $\|\nabla_g r\|^2\leq Cr^2$. To estimate $|\Delta_{g(t)} r|$ and
$\|\nabla_{g(t)} r\|_{g(t)}^2$, we recall again
$$
\fr{\partial\Gamma_{ij}^k}{\partial t}=-g^{kl}(R_{li,j}+R_{lj,i}-R_{ij,l})
$$
and thus  calculate
\begin{equation}\label{evolving laplacian}
\aligned
\fr{\partial}{\partial t}(\Delta r) & =\fr{\partial}{\partial t}(g^{ij}(\nabla^2 r)_{ij})\\ & =2g^{ki}g^{lj}h_{kl}(\nabla^2 r)_{ij} + g^{ij}g^{kl}(R_{li,j} + R_{lj,i} - R_{ij,l})\nabla_k r\\
& = 2g^{ki}g^{lj}h_{kl}(\nabla^2 r)_{ij}. \endaligned
\end{equation}
Form the fact that $C^{-1}g\leq g(t)\leq Cg$ and the
curvature estimates in Lemma \ref{curvature estimate}, we get
the desired estimates
$$|\Delta_{g(t)} r|\leq Cr\quad \mbox{and}\quad \|\nabla_{g(t)} r\|_{g(t)}^2\leq Cr^2.
$$
We consider $\bar h = r^{-\gamma} h$, $\bar{\nabla h} =
r^{-\gamma}\nabla h$, $\bar{\nabla^2 h} = r^{-\gamma}\nabla^2 h$,
and $\bar{\nabla^3 h} = r^{-\gamma} \nabla^3 h$ and calculate
\begin{eqnarray*}
\frac{\partial}{\partial t}\|\bar h\|^2 &\leq &\Delta
\|\bar h\|^2 -\|\bar{\nabla h}\|^2+ C \|\bar h\|^2,\\
\frac{\partial}{\partial t} \|\bar{\nabla h}\|^2 &\leq& \Delta \|\bar{\nabla h}\|^2 - \|\bar{\nabla^2 h} \|^2 + C (\|\bar h\|^2+\|\bar{\nabla h}\|^2),\\
\frac{\partial}{\partial t}(t\|\bar{\nabla^2 h}\|^2) &\leq & \Delta(t \|\bar{\nabla^2 h}\|^2) - t \|\bar{\nabla^3 h}\|^2\\
&& + (1 +Ct)\|\bar{\nabla^2 h}\|^2+ C (\|\bar h\|^2+\|\bar{\nabla h}\|^2).
\end{eqnarray*}
Set
$$
\varphi_1=\|\bar h\|^2+\|\bar{\nabla h}\|^2
$$
and
$$
\varphi_2=\|\bar h\|^2+\|\bar{\nabla h}\|^2+t \|\bar{\nabla^2 h}\|^2,
$$
and calculate that
$$
\frac{\partial}{\partial t} \varphi_1\leq \Delta\varphi_1+  C \varphi_1 \text{ and }
\frac{\partial}{\partial t} \varphi_2\leq \Delta\varphi_2+ C \varphi_2.
$$
Therefore
$$
\frac{\partial}{\partial t} (\re^{-Ct}\varphi_1) \leq (\re^{-Ct}\Delta\varphi_1) \text{ and }
\frac{\partial}{\partial t} (\re^{-Ct} \varphi_2) \leq (\re^{-Ct}\Delta\varphi_2).
$$
Thus the lemma follows from Lemma \ref{maximum principle}. The assumption $(iii)$ is satisfied because, for an asymptotically hyperbolic manifold,
$$
\int_{\cM} \exp{(-\alpha d(x, x_0))} dv_g < \infty
$$
for any $\alpha > n-1$ and that
$$
\frac 1C \re^{-d(x, x_0)}  \leq r \leq C \re^{- d(x, x_0)}
$$
for some constant $C$.
\end{proof}

To get the decay estimates independent of the time, we need to be
a bit more careful in the above calculation, at least, near the
boundary at the infinity. On an asymptotically hyperbolic manifold
$(\cM^n, \ g_0)$, in the coordinates at the infinity introduced by a
fixed geodesic defining function $r$ for $0 < r < r_0$, we know
that 
$$
 |\nabla_g r|_g^2  = r^2 \text{ and } \Delta_g r  = (2-n)r +
O(r^2),
$$
in the light of the normal form \eqref{normal form}. To better
estimate $\|\nabla_{g(t)}r\|^2_{g(t)}$ and $\Delta_{g(t)}r$ along the flow $g(t)$, we
calculate as follows:

\begin{lemm}\label{estimate of r along time}
Suppose that $(\cM, \ g_0)$ is an asymptoticlly hyperbolic manifold with
a fixed geodesic defining function $r$ for $0<r<r_0$. Let $g(t)$
,$t\in[0,T]$, be NRF starting from the metric $g_0$ such that
$$
\int_0^t\|h\|(\cdot, s)ds \leq \delta
$$ 
for all $t\in[0,T]$, where the number $\delta$ is independent of $T$. Then, for $0<r<r_0$ and some $C$ independet of $T$, we have
\begin{equation}\label{r-precise}
\|\nabla_{g(t)}r\|^2_{g(t)}=r^2 + C\delta r^2\quad \mbox{and}\quad \Delta_{g(t)}r=(2-n)r+C\delta r + C r^2
\end{equation}
for some $C$ independent of $T$.
\end{lemm}
\begin{proof} First, using \eqref{NRF}, we calculate
$$
\fr{\partial}{\partial t} g^{ij}\nabla_i r\nabla_j r = - 2g^{ik}h_{kl}g^{lj}\nabla_i r\nabla_j r.
$$
Hence, by the assumption, it is easy to get $\|\nabla_{g(t)}r\|^2_{g(t)} = (1+C\delta)r^2.$  Similarly, from \eqref{evolving laplacian}, 
$$
\fr{\partial}{\partial t}(\Delta r)= 2g^{ki}g^{lj}h_{kl}(\nabla^2 r)_{ij} 
$$ 
and the assumptions, it is easy to check that $\Delta_{g(t)}r=(2-n)r + C\delta r+ Cr^2.$
\end{proof}

Let us revise the calculations in the proof of Lemma \ref{t-dependent} and be more careful this time with the help of Lemma \ref{r-precise}. First let us recall the evolution equation (\ref{evolving equation 2})
$$
\frac{\partial}{\partial t}\|h\|^2=\Delta \|h\|^2-2\|\nabla
h\|^2+ 4R_{ijkl}h^{ik}h^{jl}.
$$
Let $\bar h= r^{-\gamma} h$. Then $\bar h$ satisfies the evolving equation:
$$
\frac{\partial}{\partial t}\|\bar h\|^2= r^{-2\gamma} \Delta \|h\|^2-2r^{-2\gamma}\|\nabla
h\|^2+ 4R_{ijkl}\bar h^{ik}\bar h^{jl}.
$$
We calculate
\begin{equation}\label{calculation-1}
\aligned
& r^{-2\gamma} \Delta \|h\|^2 = \Delta \|\bar h\|^2 - (\Delta r^{-2\gamma} )\|h\|^2 - 2\nabla r^{-2\gamma}\cdot\nabla \|h\|^2 \\
 = \Delta & \|\bar h\|^2 - (2\gamma(2\gamma + (n-1) ) + C\delta + C r)\|\bar h\|^2 + 4\gamma \frac {\nabla r}r \cdot r^{-2\gamma}\nabla\|h\|^2
\endaligned
\end{equation}
where
$$
\Delta r^{-2\gamma} = (2\gamma (2\gamma + (n-1)) + C \delta  + C r) r^{-2\gamma}
$$
due to Lemma \ref{estimate of r along time}. For the last term
in \eqref{calculation-1} we split into two equals and treat them
differently
$$
\aligned 4\gamma \frac {\nabla r}r  & \cdot
r^{-2\gamma}\nabla\|h\|^2 =  4\gamma \frac {\nabla r}r r^{-\gamma}
h^{ij}\cdot r^{-\gamma} \nabla h_{ij} + 2\gamma \frac {\nabla r}r
\cdot \nabla \|\bar h\|^2 + 4\gamma^2 \|\bar h\|^2\cdot\|\fr{\nabla r}{r}\|^2\\
& \leq 2\gamma^2 \|\bar h\|^2\cdot\|\fr{\nabla r}{r}\|^2 +
2r^{-2\gamma}\|\nabla h\|^2  + 2\gamma \frac {\nabla r}r \cdot
\nabla \|\bar h\|^2 + 4\gamma^2 \|\bar h\|^2\cdot\|\fr{\nabla
r}{r}\|^2
\endaligned
$$
Therefore we arrive at
\begin{equation}\label{calculation-2}
\aligned \frac {\partial}{\partial t}\|\bar h\|^2 & \leq \Delta
\|\bar h\|^2  +2\gamma \frac {\nabla r}r\cdot\nabla\|\bar h\|^2 \\ +
& 2(\gamma^2 - (n-1)\gamma + C r +C \delta) \|\bar h\|^2
+ 4R_{ijkl}\bar h^{ik}\bar h^{jl}. \endaligned
\end{equation}
Now we are ready for the main theorem of this section.

\begin{theo}\label{long time existence-asym}
Given any $n\geq 5$ and positive constants $k_0, k_1, v_0, \lambda_0$, let
$$
\gamma \in( \frac {n-1}2  -  \min\{\sqrt \lambda_0, \sqrt{\frac{(n-1)^2}4 - 2}\}, \ \frac {n-1}2 + \sqrt{\frac{(n-1)^2}4 - 2}).
$$
Suppose that $(\cM^n, \  g_0)$ is an asymptotically hyperbolic manifold of $C^2$ regularity, where $g_0$ satisfies the conditions $B(k_0, v_0,\lambda_0)$ and $\nabla Rm|\leq k_1$.
Then there exists a positive number $\delta>0$ depending only on
$n, k_0, k_1, v_0, \lambda_0, \delta$ such that NRF on $\cM^n$ starting from $g_0$ exists for all the time  if the initial metric $g_0$ is $\delta$-Einstein of order $\gamma$ and $\|\nabla h\|\leq Cr^\gamma$.
Moreover NRF converges exponentially to a
non-degenerate Einstein metric $g_\infty$ in $C^\infty$ as $t\to\infty$, and $g_\infty$ is $C^2$-conformally compact with the same conformal infinity as the initial metric $g_0$ if $\gamma > 2$.
\end{theo}

\begin{proof} In the light of Lemma \ref{lemm:analytic tool} and Theorem \ref{long time existence-pre-asym} we know NRF exists all the time and converges exponentially to a non-degenerate Einstein metric $g_\infty$ in $C^\infty$ under the assumptions of this theorem. Because $\alpha$ in Theorem \ref{long time existence-pre-asym} is any number greater than $n-1$ for the given asymptotically hyperbolic metric $g$.

To show that $g_\infty$ is conformally compact of some regularity we will apply Lemma \ref{maximum principle} to get the curvature $h$ decay both in time and in space approaching the infinity. For that we consider $\tilde h = \re^{\lambda_1 t} r^{-\gamma} h$ and calculate, from \eqref{calculation-2},
\begin{equation}\label{tilde-calculation}
\aligned
\frac d{dt}\|\tilde h\|^2 & \leq \Delta \|\tilde h\|^2 +2\gamma \frac {\nabla r}r\cdot\nabla\|\tilde h\|^2 \\ + 2( & \gamma^2 - (n-1)\gamma + C\delta + Cr) \|\tilde h\|^2 + 4R_{ijkl}\tilde h^{ik}\tilde h^{jl} + 2\lambda_1 \|\tilde h\|^2.\endaligned
\end{equation}
From the proof of Theorem \ref{long time existence-pre-asym} we also know that
$$
\|R_{ijkl} (g(t)) - R_{ijkl}(g_0)\| \leq C\delta
$$
for all $t$. Therefore we have, for $r\in (0, r_1)$ for some $r_1$,
$$
\frac d{dt}\|\tilde h\|^2  \leq \Delta \|\tilde h\|^2 +2\delta \frac {\nabla r}r\cdot\nabla\|\tilde h\|^2+ 2A\|\tilde h\|^2,
$$
where
\begin{equation}\label{the decay range}
A = \gamma^2 - (n-1)\gamma +2 + C\delta + C r + \lambda_1  \leq 0
\end{equation}
when $\delta$, $r$, and $\lambda_1$ is sufficiently small for a given
$$
\gamma \in( \frac {n-1}2  -  \min\{\sqrt \lambda_0, \sqrt{\frac{(n-1)^2}4 - 2}\}, \ \frac {n-1}2 + \sqrt{\frac{(n-1)^2}4 - 2}).
$$
To apply Lemma \ref{maximum principle} we verify the bounds at the initial time and on the boundary at $r= r_1$. To make sure $\|\tilde h\|^2$ is bounded at $r=r_1$ for all time we need to set $\lambda_1$ less than the index of exponential decay we observed from Lemma \ref{long time-asym}. Thus from Lemma \ref{maximum principle} we have
\begin{equation}\label{decay-t-s}
\|h\| (\cdot, t) \leq C \re^{-\lambda_1 t} r^{\gamma}
\end{equation}
for some positive number $\lambda_1$ and for all the time, which implies that
$$
\|r^2 g(t) - r^2 g\|_{\bar g} \leq C r^\gamma
$$
for all the time. So it is clear from here that $g_\infty$ is conformally compact with the same conformal infinity as the initial metric $g$.  For the regularity we consider
the Fermi coordinate $(r, \theta^1, \cdots, \theta^{n-1})$ near the boundary at the infinity of $\cM$ and calculate
$$
\partial_k  \bar g_{ij}(\cdot, t) - \partial_k \bar g_{ij}(\cdot, 0) = -4 r\partial_k r \int_0^t h_{ij} - 2 r \int_0^t r\partial_k h_{ij}
$$
where $\partial_k = \frac \partial{\partial \theta^k}$, $r=
\theta^0$, and $\bar g (t) = r^2 g(t)$. By direct calculation, one
easily gets that $\nabla_k h_{ij}=\partial_k h_{ij}+ C
r^{-3}\|h\|.$Hence
$$
\|\partial_k \bar g_{ij}(\cdot, t) - \partial_k\bar g_{ij}(\cdot, 0)\|_{\bar g} \leq Cr^{\gamma -1} + C r^{-1}\int_0^t \|\nabla h\|.
$$
Therefore to say $\bar g(t)$ is $C^1$ it suffices to get the estimates
\begin{equation}\label{derivatives-h}
\|\nabla h\| (\cdot, t)\leq C \re^{-\lambda_1 t}r^\gamma
\end{equation}
for all the time $t$ for some $\gamma > 1$. Moreover to show $g(t)$ is $C^2$ conformally compact we need, in addition, 
\begin{equation}\label{2-derivatives-h}
\|\nabla^2 h\| (\cdot, t)\leq C \re^{-\lambda_1 t}r^\gamma
\end{equation}
for all $t$ and some $\gamma > 2$. From Lemma \ref{lemma3} we know \eqref{derivatives-h} and \eqref{2-derivatives-h} hold for $t$ greater than some small number. Therefore the issue is to get decay of $\|\nabla h\|$ and $\|\nabla^2 h\|$ in space for small $t>0$. In a similar way it is easily seen that to show that $\bar g(t)$ is of $C^2$ it suffices to show that
$$
\|\nabla h\| (\cdot, t) \leq C r^\gamma \text{ and } \|\nabla^2 h\|(\cdot, t) \leq \frac C{\sqrt t} r^\gamma
$$
for  $t > 0$ small, when $\gamma > 2$.  The proof is complete because of  Lemma \ref{t-dependent}.
\end{proof}

\begin{rema} The reason that we are only concerned with $C^2$ regularity for the conformal compactness of the metric $g_\infty$ in the above theorem is 
because the regularity theorem in \cite{CDLS}, which  states
that a conformally compact Einstein metric is smooth in all even
dimensions and polyhomogeneous in odd dimensions greater than $3$ if
it is conformally compact of $C^2$ regularity and the conformal
infinity $(\partial\cM, \ [\hat g])$ is smooth.
\end{rema}


\section{Perturbation existence of CCE metrics}

In this section we would like to apply Theorem \ref{long time existence-asym} to recapture some works in \cite{GL} \cite{Lee} \cite{Bi}. Our approach is to construct a good initial metric nearby a given non-degenerate conformally compact Einstein metric and to apply Theorem \ref{long time existence-asym}. The construction is based on the metric expansions of Fefferman and Graham.
We glue together the given metric away from the boundary and the metric on collar
neighborhood of the boundary from the expansion of Fefferman and Graham. The issues are to make sure the glued metrics satisfy the assumptions on the initial metric in Theorem \ref{long time existence-asym}.

Suppose that $(\cM^n, \ g)$ is a conformally compact Einstein manifold with the conformal infinity $(\partial\cM, \ [\hat g])$. Suppose that $r$ is the geodesic defining function associated with the conformal representative $\hat g\in [\hat g]$ on $\partial\cM$.
Then the metric expansion is given as follows (cf. \cite{FG}):
$$
\aligned
g_r & = \hat g + g^{(2)}r^2 + \cdots +  g^{(n-3)}r^{n-3} + hr^{n-1} \log r + g^{(n-1)} r^{n-1} + \cdots \\
& = \hat g + g^{(2)}r^2 + \cdots +  g^{(k)}r^{k} + t^{(k)}[g],
\endaligned
$$
for $0 \leq k\leq n-3$, when $n-1$ is even,
$$
\aligned
g_r & = \hat g + g^{(2)}r^2 + \cdots + g^{(n-2)}r^{n-2} + g^{(n-1)} r^{n-1} + \cdots\\
& = \hat g + g^{(2)}r^2 + \cdots +  g^{(k)}r^{k} + t^{(k)}[g],
\endaligned
$$
for $0 \leq k \leq n-2$, when $n-1$ is odd,
where
\begin{itemize}
\item
$g^{(2i)} $ for $2i< n-1$ are local invariants of $(\partial\cM^{n-1}, \hat g)$;
\item
$h$ and $tr g^{(n-1)}$ ($n-1$ even) are also local invariant of
$(\partial\cM^{n-1}, \hat g)$;
\item
$h$ and $g^{(n-1)}$ ($n-1$ odd) are trace free;
\item
$g^{(n-1)}$ ($n-1$ odd) and trace-free part of $g^{(n-1)}$ ($n-1$ even)
are nonlocal.
\end{itemize}
For instance,  
$$
g^{(2)} = -\frac 1{n-3}(\hat Ric - \frac {\hat R}{2(n-2)}\hat g).
$$
To construct a candidate to be the right initial metric to apply Theorem \ref{long time existence-asym}, whose conformal infinity is a perturbation of that of a given conformally compact Einstein metric $g$, we set
\begin{equation}\label{full expansion}
 g^{k,\nu}_r = \hat g_\nu + g^{(2)}_\nu r^2 + \cdots +  g^{(k)}_\nu r^{k}  + t^{(k)}[g]
\end{equation}
where $\hat g_\nu$ is a perturbation of $\hat g$,  and $g^{(2i)}_\nu= g^{(2i)}[\hat g_\nu]$, $2i\leq k$, are corresponding curvature terms of $\hat g_\nu$ as given in the metric expansion in \cite{FG}. Next let $\phi$ be a cut-off function of the variable $r$ such that $\phi =0$ when $r \geq \nu_2$ and $\phi = 1$ when $r \leq \nu_1$, where $\nu_1 < \nu_2$ are chosen later. We therefore have the candidate
\begin{equation}\label{candidate}
g_{k, \nu}^\phi = r^{-2} (dr^2 + (1 - \phi) g_r +  \phi g^{k,\nu}_r).
\end{equation}

Immediately we see that
$$
\|g_{k,\nu}^\phi - g\|_g \leq C \|\hat g_\nu - \hat g\|_{C^{k}}.
$$
For the convenience of readers, we include all curvature calculations in the appendix. From the calculations in the appendix we observe that
$$
\|\Gamma^l_{ij}[g_{k,\nu}^\phi] - \Gamma^l_{ij}[g]\|_g \leq C  \|\hat g_\nu - \hat g\|_{C^{k+1}}
$$
$$
\| Rm[g_{k,\nu}^\phi] -  Rm[g]\|_g \leq C \| \hat g_\nu - \hat g\|_{C^{k+2}}
$$
and
\begin{equation}
\|\nabla Rm[g_{k, \nu}^\phi] - \nabla Rm[g]\| \leq  C\|\hat g_\nu-
\hat g\|_{C^{k+2}}.
\end{equation}
Hence the metric $g_{k,\nu}^\phi$ satisfies the conditions $B(k_0, v_0, \lambda_0)$ and $|\nabla Rm|\leq k_1$ for some positive numbers $k_0, k_1, v_0$, and $\lambda_0$,  if the conformally compact Einstein metric $g$ is non-degenerate and $\|\hat g^\nu - \hat g\|_{C^{k+2}}$ is sufficiently small.

It is also immediate from \eqref{candidate} that
$$
h [g_{k, \nu}^\phi] = h[g] = 0 \text{ in $\{r > \nu_2\}$}.
$$
It is important to observe from \eqref{h-curvature} that, for $\nu_1 < r < \nu_2$,
$$
\|h\|[g_{k, \nu}^\phi] \leq C r^2 \|\hat g_\nu - \hat g\|_{C^{k+2}},
$$
and, for $r<\nu_1$,
$$
\|h\| [g_{k, \nu}^\phi] = Cr^{k+2} \|\hat g_\nu - \hat g\|_{C^{k+2}}
$$
and
\begin{equation}
\|\nabla h\|[g_{k, \nu}^\phi] = C r^{k+2} \|\hat g_\nu - \hat g\|_{C^{k+2}}.
\end{equation}
Thus, as a consequence of Theorem \ref{long time existence-asym}, we are able to recover some of the works in \cite{Lee} \cite{Bi} as follows:

\begin{theo}\label{general perturbation} Let $(\cM^n, \ g)$  be a conformally compact Einstein manifold of regularity $C^2$ with a 
smooth conformal infinity $(\partial\cM, \ [\hat g])$.  Assume that $g$ is of the non-degeneracy $\lambda_0$ as defined in Definition \ref{def.ofconditionB}.  Suppose that
\begin{equation}\label{eq:lambda}
\max\{2, \frac {n-1}2 - \sqrt{\lambda_0}\} < k+2.
\end{equation}
for some even $k \leq n-1$. Then, if a smooth metric $[\hat g_\nu]$ is a sufficiently small $C^{k+2}$-perturbation of $[\hat g]$, then there is a $C^2$-conformally compact Einstein metric on $\cM$ whose conformal infinity is $[\hat g_\nu]$.
\end{theo}

\begin{proof} First of all, from the above discussion, it is clear that $\hat g_{k,\nu}^{\phi}$ satisfies the conditions $B(k_0, v_0, \lambda_0)$ and $|\nabla Rm|\leq k_1$, where the constants $k_0, k_1, v_0, \lambda_0$ are close to those of the metric $g$, provided that $\hat g_\nu$ is a smooth sufficiently small $C^{k+2}$-perturbation of $\hat g$. We claim that, for any $2<\gamma < k+2$,
$$
\|h\|[g_{k,\nu}^\phi] \leq \epsilon r^\gamma \text{ and } \|\nabla h\|[g_{k, \nu}^\phi] \leq Cr^\gamma,
$$
provided that $\hat g_\nu$ is a smooth sufficiently small $C^{k+2}$-perturbation of $\hat g$. In fact, for an appropriately small $\nu_1 >0$, for any $\epsilon>0$,
$$
\|h\|[g_{k, \nu}^\phi] \leq \epsilon r^\gamma,
$$
for all $r< \nu_1$, provided that $\hat g_\nu$ is a smooth sufficiently small $C^{k+2}$-perturbation of $\hat g$. We also know that
$$
\|h\|[g_{k,\nu}^\phi] \leq C r^2\|\hat g_\nu - \hat g\|_{C^{k+2}} \leq \epsilon \nu_1^\gamma
$$
for $r \geq \nu_1$, if $\|\hat g_\nu - \hat g\|_{C^{k+2}}$ is small enough. We will take $\nu_2 = 2\nu_1$. Notice that $\|\nabla \phi\|_g \leq C$. Therefore we can have
\begin{equation}\label{main-condition}
\|h\|[g_{k, \nu}^\phi] \leq \epsilon r^\gamma \text{ and } \|\nabla h\|[g_{k, \nu}^\phi] \leq \epsilon r^\gamma \text{ on $\cM$}
\end{equation}
when $\|\hat g_\nu - \hat g\|_{C^{k+2}}$ is sufficiently small. Of course we will use the number $\gamma$ such that 
$$
\gamma < \frac {n-1}2 + \sqrt{\frac {(n-1)^2}4 - 2}.
$$
Thus Theorem \ref{long time existence-asym} applies and gives the existence of $C^2$-conformally compact Einstein metrc on $\cM$ with the conformal infinity $(\partial\cM, \ [\hat g_\nu])$.
\end{proof}

We would like to remark that the reason we need to consider $C^{k+2}$ perturbations is because of \eqref{delta} in Lemma \ref{lemma:pinch with order}. With a little more careful calculations and the use of H\"{o}lder norm we can obtain the following: 

\begin{theo}\label{two plus alpha}
Let $(\cM^n, \ g)$, $n\geq 5$,   be a conformally compact Einstein manifold of regularity $C^2$ with a 
smooth conformal infinity $(\partial\cM, \ [\hat g])$.  And suppose that the non-degeneracy
$$
\sqrt \lambda>\frac {n-1}2 -2
$$
for $(\cM^n, \ g)$. Then, for any smooth metric $\hat g_\nu$ on $\partial\cM$, which is sufficiently
$C^{2,\alpha}$ close to  some $\hat g\in [\hat g]$ for any $\alpha\in(0,1)$, there
is a conformally compact Einstein metric $g_{\nu}$ on $\cM$ which can be $C^2$
conformally compactfied with the conformal infinity  $[\hat g_\nu]$.
\end{theo}

To take advantage of the H\"{o}lder continuity in the expansions we use the following simple analytic lemma.
For simplicity we will state in the context of a piece of flat boundary in Euclidean space, but it is easily seen to hold 
in our context when the coordinate $x_n$ is replaced by the geodesic defining function of asymptotically hyperbolic
manifolds.

\begin{lemm}\label{bdrregularity}
Let $B(2R) =\{x\in \mathbf{R}^n:|x| \leq 2R, \ x_n=0\}$ and $\phi\in C^{0,\alpha}(B(2R))$. Then we can extend $\phi$ to  a smooth function $u$ 
on $B(R)\times(0, R)\subset \mathbf{R}^n$ such that
$$
|u(x',x_n)-\phi(x')|=C x_n^{\alpha},
$$
and
$$
|D^\beta u(x', x_n)|= Cx_n^{-m+\alpha},
$$
where $x=(x',x_n)$,  $|\beta|=|(\beta' , \beta_n)| = m$, and $m \geq 1$. 
\end{lemm}

\begin{proof} We will use a mollifier to
extend the boundary function and show that the smoothened function
satisfies the desired estimates. Precisely, let $\tau\in
C^{\infty}_0(R^{n-1}), \tau \geq 0$ be a mollifier with  support in $B(1)\subset R^{n-1}$ and satisfy
$$
\int_{R^{n-1}}\tau =1.
$$
We then define
$$
u(x',x_n)=x_n^{1-n}\int_{\omega}\tau(\frac{x'-y'}{x_n})\phi(y')dy' \in C^\infty(B(R)\times(0, R)).
$$
Hence it is easily calculated that
\begin{eqnarray*}
u(x',x_n)-\phi(x')&=&x_n^{1-n}
\int_{\omega}\tau(\frac{x'-y'}{x_n})(\phi(y')-\phi(x'))dy'\\
&=& \int_{B(1)} \tau(z')(\phi(x'- x_n z')-\phi(x'))dz'\nonumber\\
&\leq& C  x_n^{\alpha}\nonumber
\end{eqnarray*}
When $m\geq 1$ we always have
\begin{eqnarray*}
D^\beta u(x',x_n) = x_n^{1 - n - m}\int P_\beta (\frac {x'-y'}{x_n})(\phi(y') - \phi(x')) dy',
\end{eqnarray*}
where $P_\beta \in C_0^{\infty}$ with support in $B(1)$.  For instance, 
$$
P_{(0, 1)} (z') = (1-n)\tau (z') - \nabla\tau (z')\cdot z'
$$
and
$$
P_{(\beta', 0)} (z') = D^{\beta'} \tau (z').
$$
Therefore we can conclude that
$$
|D^\beta u(x', x_n)| \leq C x_n^{-m +\alpha}.
$$
\end{proof}

Now we can use the above lemma to construct a good initial metric near by a given non-degenerate conformally compact Einstein metric
in $C^{2, \alpha}$ topology.

\begin{lemm}\label{approximate AHE metric}
Let $(\cM^n, \ g)$  be a conformally compact Einstein manifold of regularity $C^2$ with a 
smooth conformal infinity $(\partial\cM, \ [\hat g])$.  Assume $r$ is the geodesic defining
function associated with $\hat g\in [\hat g]$. Let $\hat{g}_{\nu}$ be a smooth perturbation of  $\hat g$ on
$\partial\cM$ such that
$$
\|\hat{g}_{\nu}-\hat g\|_{C^{2,\alpha}}\leq \varepsilon.
$$
Then there is a smooth AH metric $g_0$ on $\cM^n$with 
conformal infinity $[\hat{g}_{\nu}]$, and
\begin{enumerate}
\item $\|g_0-g\|_g\leq C\varepsilon$ and $\|\Gamma [g_0] - \Gamma [g]\| \leq C \varepsilon$ ;
\item $\|h(g_0)\|_g=C r^{2+\alpha}$ and $\|\nabla h(g_0)\|_g=C r^{2+\alpha}$;
\item $\|Rm [g_0] - Rm[g]\|_g\leq C\varepsilon$ and $\|\nabla Rm [g_0] - \nabla Rm[g]\|_g \leq C$,
\end{enumerate}
near the boundary, where $C$ is a constant independent of $g_0$.
\end{lemm}

\begin{proof}
First let 
$$
\hat T_\nu = -\frac1{n-3}(Ric(\hat{g}_{\nu})-\frac{R(\hat{g}_{\nu})}{2(n-2)}\hat{g}_{\nu})
$$
the Schouten tensor for $\hat g_\nu$, which is close to the Schouten tensor $\hat T$ of the metric $\hat g$ in $C^{0, \alpha}$ on the boundary
$\partial\cM$. We apply Lemma \ref{bdrregularity} to extend $\hat T_\nu$
to $T$,  which is defined in a neighborhood of $\partial\cM$. Then we construct
$$
g_0= r^{-2}(dr^2+ \hat{g}_{\nu} + r^2 T)
$$
near the boundary, and defined to be $g$ inside of $\cM$ as we did in the proof Theorem \ref{general perturbation}.
Therefore the lemma follows from Lemma \ref{bdrregularity} and calculations in the appendix.
\end{proof}

Finally the proof of Theorem \ref{two plus alpha} goes with little change from the proof of Theorem \ref{general perturbation} after Lemma \ref{approximate AHE metric}. 
It is good to realize that Theorem \ref{two plus alpha} recovers the perturbation theorem of Graham and Lee in \cite{GL}, because a perturbation of the standard hyperbolic
metric on the unit Euclidean ball $B^n$ as given in Theorem \ref{two plus alpha} clearly has the non-degeneracy close to that of hyperbolic metric, which is $\frac {(n-1)^2}4$.


\appendix

\section{Curvature for AH metrics}

Suppose that $(\cM, \ g)$ is an asymptotically hyperbolic manifold and that $r$ is the geodesic defining function associated with the representative $\hat g$ of the conformal infinity $(\partial\cM, \ [\hat g])$. Hence in the Fermi coordinate at the infinity boundary we have the metric in norm form
\begin{equation}
g = r^{-2} \bar g = r^{-2}(dr^2 + g_r).
\end{equation}
We will use the convention that the Latin letters stand for the index: $1, 2, \cdots, n$ and the Greek letters stand for the index: $1, 2, \cdots, n-1$. We will identify $r = x^n$. Recall that the Riemann curvature tensor is given as
\begin{equation}
\aligned
R_{ijkl} & = \frac 12(-\partial_j\partial_lg_{ik} - \partial_i\partial_kg_{jl} + \partial_i\partial_l g_{jk} +\partial_j\partial_k g_{il})\\ & - g^{mn}([ik, m][jl, n] - [il, m][jk, n])\endaligned
\end{equation}
and the Ricci curvature tensor is given as
\begin{equation}
R_{ik} = g^{jl}R_{ijkl},
\end{equation}
where the Chrsitoffel symbol of second kind is given as
\begin{equation}
[ij,k] = \frac 12(\partial_ig_{jk} + \partial_j g_{ik} - \partial_k g_{ij}).
\end{equation}
We start with computing Christofell symbols of second kind.
\begin{equation}
[nn,n] =  -\frac 1{r^3}
\end{equation}
\begin{equation}
[nn, \alpha] = [n\alpha,n] = 0
\end{equation}
\begin{equation}
[\alpha\beta, n] = - [n \alpha, \beta] = r^{-3}\bar g_{\alpha\beta} - \frac 12r^{-2}\bar g_{\alpha\beta}'
\end{equation}
and
\begin{equation}
[\alpha\beta,\gamma] = r^{-2}[\alpha\beta,\gamma][\bar g].
\end{equation}
Then we calculate the Riemann curvature tensor.
$$
R_{nnnn} = R_{nnnk} = R_{nnij} = 0
$$
\begin{equation}
R_{n\alpha n\beta}  = -r^{-4} \bar g_{\alpha\beta}+\frac 12 r^{-3}\bar g_{\alpha\beta}' + \frac 14 r^{-2}\bar g^{\gamma\delta}\bar g_{\alpha\delta}'\bar g_{\gamma\beta}' -\frac 12r^{-2}\bar g_{\alpha\beta}''
\end{equation}
\begin{equation}
R_{n\gamma\alpha\delta} = \frac 12 r^{-2}( -\bar \nabla_\alpha \bar g_{\gamma\delta}' + \bar\nabla_\delta \bar g_{\alpha\gamma}')
\end{equation}
and
\begin{equation}
\aligned
R_{\alpha\gamma\beta\delta}  & = - r^{-4}(\bar g_{\alpha\beta}\bar g_{\gamma\delta} - \bar g_{\alpha\delta}\bar g_{\gamma\beta})  - \frac 14 r^{-2} (\bar g_{\alpha\beta}'\bar g_{\gamma\delta}' - \bar g_{\alpha\delta}'\bar g_{\gamma\beta}')\\+  & \frac 12 r^{-3}(\bar g_{\alpha\beta}\bar g_{\gamma\delta}' + \bar g_{\gamma\delta}\bar g_{\alpha\beta}' - \bar g_{\alpha\delta}\bar g_{\gamma\beta}' - \bar g_{\gamma\beta}\bar g_{\alpha\delta}') +  r^{-2}\bar R_{\alpha\gamma\beta\delta}
\endaligned
\end{equation}
Therefore we can get the Ricci curvature.
$$
R_{nn} =  - (n-1) r^{-2} +\frac 12 r^{-1}\bar g^{\alpha\beta}\bar g_{\alpha\beta}' +\frac 14 \bar g^{\alpha\beta}\bar g^{\gamma\delta}\bar g_{\alpha\delta}'\bar g_{\gamma\beta}' - \frac 12 \bar g^{\alpha\beta}\bar g_{\alpha\beta}''.
$$
$$
R_{n \alpha}  = \frac 12 \bar g^{\gamma\delta}( -\bar \nabla_\alpha \bar g_{\gamma\delta}' + \bar\nabla_\delta \bar g_{\alpha\gamma}')
$$
and
$$
\aligned R_{\alpha\beta}  & = -(n-1)r^{-2}\bar g_{\alpha\beta} -
\frac 12 \bar g_{\alpha\beta}''+(\frac n2-1)r^{-1}\bar
g_{\alpha\beta}'+ \frac 12 r^{-1}\bar g^{\gamma\delta}\bar
g_{\gamma\delta}' g_{\alpha\beta}\\ & \quad  - \frac 14 \bar
g^{\gamma\delta}\bar g_{\gamma\delta}' \bar g_{\alpha\beta}' +\frac
12 \bar g^{\gamma\delta}\bar g_{\alpha\delta}'\bar g_{\beta\gamma}'
+ Ric_{\alpha\beta}(g_r).
\endaligned
$$
To summarize, for the curvature $h = Ric + (n-1)g$, we have
\begin{equation}\label{h-curvature}
\aligned
h_{nn} &  = \frac 12 r^{-1}\bar g^{\alpha\beta}\bar g_{\alpha\beta}' + \frac 14 \bar g^{\alpha\beta}\bar g^{\gamma\delta}\bar g_{\alpha\delta}'\bar g_{\gamma\beta}' -\frac 12 \bar g^{\alpha\beta}\bar g_{\alpha\beta}'' \\
h_{n\alpha} & = \frac 12 \bar g^{\gamma\delta}( - \bar\nabla_\alpha \bar g_{\gamma\delta}' + \bar\nabla_\delta \bar g_{\alpha\gamma}')\\
h_{\alpha\beta} & = - \frac 12 \bar g_{\alpha\beta}''+(\frac
{n}2-1)r^{-1}\bar g_{\alpha\beta}'+ \frac 12 r^{-1}\bar
g^{\gamma\delta}\bar g_{\gamma\delta}' g_{\alpha\beta}  \\ &
\quad\quad\quad -\frac 14 \bar g^{\gamma\delta}\bar
g_{\gamma\delta}' \bar g_{\alpha\beta}' +\frac 12 \bar
g^{\gamma\delta}\bar g_{\alpha\delta}'\bar g_{\beta\gamma}' +
Ric_{\alpha\beta}(g_r).
\endaligned
\end{equation}

\section{Evolution equations for curvatures}

We collect some basic evolution equations of geometric quantities under NRF in the following lemma for the convenience of the readers.  

\begin{lemm}\label{basicequations}
Along the normalized Ricci flow (\ref{NRF}), we have:

\begin{eqnarray}\label{evolving equation 1} \frac{\partial}{\partial
t}h_{ij} &=&\Delta
h_{ij}+2R_{ipjq}h^{pq}-2h_{ip}h_j^p=-(\Delta_L+2(n-1))h_{ij}\\
\label{evolving equation 2} \frac{\partial}{\partial
t}\|h\|^2&=&\Delta \|h\|^2-2\|\nabla
h\|^2+ 4R_{ipjq}h^{pq}h^{ij}\\
\label{evolvingequation 3} \frac{\partial}{\partial
t}R_{ijk}^{l}&=&-g^{lp}\{\nabla_{i}\nabla_{j}h_{kp}+\nabla_{i}\nabla_{k}h_{jp}
-\nabla_{i}\nabla_{p}h_{jk}\\
&&-\nabla_{j}\nabla_{i}h_{kp}-\nabla_{j}\nabla_{k}h_{ip}+\nabla_{j}\nabla_{p}h_{ik}
\}\nonumber\\
\label{evolving equation 4} \frac{\pr}{\pr t}\|
 Rm\|^2 &=& \Delta \|
 Rm\|^2 -2\|\nabla
 Rm\|^2 +Rm\ast Rm  + h\ast Rm \ast Rm \\
 \label{evolving equation 5} \frac{\pr}{\pr t}\|\nabla
Rm\|^2 &=& \Delta \|\nabla Rm\|^2 - 2\|\nabla^2  Rm\|^2 +Rm\ast\nabla Rm\ast\nabla Rm \\
 \label{evolvingequation 6} \frac{\pr}{\pr t}\|\nabla^2
 Rm\|^2 &=& \Delta \|\nabla^2
 Rm\|^2 -2\|\nabla^3
 Rm\|^2 +Rm\ast\nabla^2 Rm\ast\nabla^2 Rm\\
 &&+\nabla Rm\ast\nabla Rm\ast\nabla^2 Rm\nonumber\\
 \label{evolvingequation 7}
  \frac{\partial}{\partial
  t}\|\nabla h\|^2&=&\Delta\|\nabla h\|^2-2\|\nabla^2 h\|^2
  +h\ast\nabla h\ast\nabla Rm\\
  &&+ Rm\ast\nabla h\ast\nabla h\nonumber\\
  \label{evolvingequation 8}
  \frac{\partial}{\partial
  t}\|\nabla^2 h\|^2&=&\Delta\|\nabla^2 h\|^2-2\|\nabla^3 h\|^2
  +\nabla Rm\ast\nabla h\ast\nabla^2 h\\
  &&+ Rm\ast\nabla^2 h\ast\nabla^2 h
  +h\ast\nabla^2 Rm\ast\nabla^2 h\nonumber
  \end{eqnarray}
  Here, $\ast$  denotes linear combinations (including contractions
with $g(t)$ and its inverse ${g(t)}^{-1}$).
  \end{lemm}

\begin{proof}[Proof of Lemma \ref{basicequations}]
For the  calculations of (\ref{evolving equation 1}) and
(\ref{evolvingequation 3}), one may refer to (2.31) and (2.66) in
\cite{CLN}. Notice that
$$\frac{\partial}{\partial t}\nabla A-\nabla\frac{\partial}{\partial t} A = A\ast\nabla h$$
and that 
$$\nabla(\Delta A)-\Delta(\nabla A)=\nabla Rm\ast A +Rm\ast \nabla
A$$  for tensor fields $A$. One then may calculate all the evolution equations in the lemma.
\end{proof}



\begin{thebibliography}{abcde}

\bibitem[Ba]{Ba} E. Bahuaud {\it Ricci flow of conformally compact metrics}, arXiv:1011.2999v1

\bibitem[Bam1]{Bam1} R. Bamler {\it Stability of hyperbolic manifolds with cusps under Ricci flow}, arXiv:1004.2058v1

\bibitem[Bam2]{Bam2} R. Bamler {\it Stability of symmetric spaces of noncompact type under Ricci flow}, arXiv:1011.4267v1

\bibitem[Besse]{Besse}  A. L. Besse, {\it Einstein manifolds}, Springer-Verlag, Berlin, (1987).

\bibitem[Bi]{Bi}  O. Biquard, {\it Einstein deformations of hyperbolic metrics}, Surveys in Differential geometry: Essays on Einstein manifolds. Int. Press,Boston, MA, 1999,
pp.235-246.

\bibitem[CDLS]{CDLS} P. Chru\'sciel, E. Delay, J.M.  Lee and D. Skinner, {\it Boundary regularity of conformally compact Einstein metrics}, J. Diff. Geom.  69 (2005), no. 1, 
111 - 136. 

\bibitem[EH]{EH} K. Ecker and G. Huisken, {\it Interior estimates for hypersurfaces moving by mean curvature}, Inventions Math., vol. 105 (1991), 547-569.

\bibitem[FG]{FG} C. Fefferman and C.R. Graham, {\it Conformal invariants}, Elie Cartan et les Math$\grave{e}$matiques
d'aujourd'hui, Asterisque, (1985), 95-116.

\bibitem[Che]{Che} J. Cheeger, {\it Finiteness theorems for Riemannian manifolds}, Am. J. Math. 92 (1970), 61-74.

\bibitem[GL]{GL} C. R. Graham and J. M.  Lee,  {\it Einstein metrics with prescribed conformal infinity on the ball}, Adv. Math.  87 (1991), 186-255.

\bibitem[Gri]{Gri} A. Grigor$'$yan, {\it Heat kernel upper bounds on complete non-compact manifolds}, Revista Math Ibeoamericana, Vol. 10 (1994), No. 2, 395-452.

\bibitem[Heb]{Heb} E. Hebey, {\it Sobolev Spaces on Riemannian manifolds}, Lecture Notes in Mathematics 1635 Springer 1996.

\bibitem[Lee]{Lee} J. M. Lee, {\it Fredholm operators and Einstein metrics on conformally compact manifolds}, Memoirs of the American Mathematical Society, A.M.S., 2006.

\bibitem[LT] {LT} G. Liao and L. Tam, {\it On the heat equation for harmonic maps from non-compact manifolds}, Pacific J. math. vol 153 (1992), no. 1 129-145.

\bibitem[LY]{LY} H. Li and H. Yin, {\it On stability of the hyperbolic space form under the normalized Ricci flow}, Int. Math. Res. Not,  (2010), no. 15, 2903 - 2924.

\bibitem[CLN]{CLN} B. Chow, P. Lu and L. Ni, {\it Hamilton's Ricci flow}, Lectures in Contemporary Mathmatics 3, AMS, (1998).

\bibitem[Shi]{Shi} W. Shi, {\it Ricci deformation of the metric on complete noncompact Riemannian manifolds}, J. Diff. Geom. , Vol. 30 (1989), No.2, 303 - 394.

\bibitem[SSS]{SSS} O.C. Schn\"urer, F. Schulze and M. Simon, {\it Stability of hyperbolic space under Ricci flow}, arXiv:1003.210.

\bibitem[Su]{Su} V. Suneeta, {\it investigating the off-shell stability of anti-de Sitter space in string theory}, Classical Quantum Gravity, vol.  26 (2009), no. 3, 035023 (18 pp).

\bibitem[Ye]{Ye} R.G. Ye, {\it Ricci flow, Einstein metrics,  and space forms}, Trans. A.M.S., Vol. 338 (1993), No.2, 871-896.

\end{thebibliography}
\end{document}